\everymath{\displaystyle}
\documentclass{article}

\usepackage{amsmath}
\usepackage{amsthm}
\usepackage{amsfonts}
\usepackage{cases}
\usepackage[per-mode=symbol,
            inter-unit-product=\cdot,
	   retain-unity-mantissa=false]{siunitx}
\usepackage{physics}
\usepackage{dsfont}
\usepackage{todonotes}
\usepackage[hypertexnames=false,
            colorlinks=true,
            pdfborderstyle={/S/U/W 0},
            citecolor=blue]{hyperref}
\usepackage[capitalize]{cleveref}
\usepackage{subcaption}
\usepackage{placeins}
\usepackage{pdflscape}
\usepackage{tabularx}
\usepackage{adjustbox}
\usepackage{multirow}
\usepackage{comment}
\usepackage{mathtools}
\usepackage{bm}
\theoremstyle{definition}
\newtheorem{hp}{Assumption}

\theoremstyle{definition}
\newtheorem{rem}{Remark}
\Crefmultiformat{rem}{Remarks~#2#1#3}{ and~#2#1#3}{, #2#1#3}{ and~#2#1#3}

\newcommand{\ale}{\text{ALE}}
\renewcommand{\k}{\text{k}}

\newcommand{\la}{\text{LA}}
\newcommand{\lv}{\text{LV}}
\newcommand{\ao}{\text{AA}}
\newcommand{\mv}{\text{MV}}
\newcommand{\av}{\text{AV}}

\newcommand{\siresistance}{\si{\mmhg\second\per\milli\litre}}
\newcommand{\sicapacitance}{\si{\milli\litre\per\mmhg}}
\newcommand{\siinductance}{\si{\mmhg\second\squared\per\milli\litre}}
\newcommand{\siflowrate}{\si{\milli\litre\per\second} }
\newcommand{\sielastance}{\si{\mmhg\per\milli\litre}}

\DeclareSIUnit\mmhg{mmHg}

\title{Modeling isovolumetric phases in cardiac flows by an Augmented Resistive Immersed Implicit Surface Method}

\author{Alberto Zingaro$^{\ddagger,1}$,
				Michele Bucelli$^{\ddagger,1}$,
        Ivan Fumagalli$^1$,\\
        Luca Dede'$^1$,
        Alfio Quarteroni$^{1, 2}$}

\date{\footnotesize \textsuperscript{1} MOX, Department of Mathematics, Politecnico di Milano, P.zza Leonardo da Vinci 32, 20133 Milan, Italy \\
      \textsuperscript{2} Mathematics Institute, EPFL, Av. Piccard, CH-1015 Lausanne, Switzerland (Professor Emeritus)}

\begin{document}

\maketitle
\makeatletter{\renewcommand*{\@makefnmark}{}\footnotetext{$^\ddagger$ Alberto Zingaro and Michele Bucelli  equally contributed to this work.}\makeatother
}

\begin{abstract}
 {
A major challenge in the computational fluid dynamics modeling of the heart function is the simulation of isovolumetric phases when the hemodynamics problem is driven by a prescribed boundary displacement.
During such phases, both atrioventricular and semilunar valves are closed: consequently, the ventricular pressure may not be uniquely defined, and spurious oscillations may arise in numerical simulations.}
In this paper, we propose a suitable modification of the Resistive Immersed Implicit Surface (RIIS) method (Fedele et al., 2017) by introducing a reaction term to correctly capture the pressure transients during isovolumetric phases. The method, that we call Augmented RIIS (ARIIS) method, extends the previously proposed ARIS method (This et al., 2020) to the case of a mesh  {which} is not body-fitted to the valves. We test the proposed method on two different benchmark problems, including a new simplified problem that retains all the characteristics of a heart cycle. We apply the ARIIS method to a  {fluid dynamics} simulation of a realistic left heart geometry,  {and we show that ARIIS allows to correctly simulate isovolumetric phases, differently from standard RIIS method.}
\end{abstract}


\section{Introduction}
During the heart cycle, there are two phases in which all cardiac valves are closed and the action  {of the ventricular displacement} affects blood pressure without a net flow.
 {In the left ventricle (the same applies for the right part of the heart), during the \emph{isovolumetric contraction}, the intraventricular pressure raises up to the point in which the aortic valve open for the systolic ejection, while in the \emph{isovolumetric relaxation} the ventricular pressure decreases until reaching the atrial one, thus leading to the opening of the mitral valve.}
Cardiac valve dynamics is mainly driven by transvalvular pressure drop \cite{formaggia2006numerical}.
Hence, an accurate modeling of the isovolumetric phases in which the intraventricular pressure undergoes rapid changes is an essential prerequisite to capture valve opening and closing, and to properly model their effect on the flow.

The behavior of blood pressure in the heart chambers is determined by the contraction and relaxation of the myocardium.
With this in mind, Fluid-Structure Interaction (FSI) models coupling the blood flow with the heart mechanics have been proposed in the literature \cite{brenneisen2021sequential,bucelli2022partitioned,cheng2005fluid,khodaei2021personalized,nordsletten2011fluid,zhang2001analysis, bucelli2022mathematical}, or even more realistic electrophysiology-mechanics-hemodynamics models as in, e.g., \cite{santiago2018fully,watanabe2002computer,viola2020fluid,choi2015new}.
However, these coupled models typically entail a high computational cost, and they require a challenging calibration of  {a huge number of} physical parameters, especially in pathological conditions.
Because of this, uncoupled (or one-way coupled) approaches have been proposed, to address the sole Computational Fluid Dynamics (CFD) component of the system, with the ventricular displacement prescribed as data coming from analytical functions \cite{tagliabue2017fluid,zingaro2021hemodynamics,domenichini2005three, baccani2002vortex,mittal2016computational,corti2022impact}, clinical measurements \cite{fumagalli2020image,this2020pipeline,chnafa2014image,masci2020proof}, or from electromechanical simulations \cite{augustin2016patient,this2020augmented,zingaro2022geometric}.
Such models mainly differ in the treatment of the valve geometry and dynamics.
Mesh-conforming approaches are based on a classical Arbitrary Lagrangian-Eulerian formulation of the flow equations \cite{chengAleFsiAv,espinoAleValve,xale,zhangAleValve,vergaraFSIAV}, and they include the Resistive Immersed Surface (RIS) method \cite{fernandez2008numerical,astorino2012robust} and different XFEM/cutFEM methods \cite{alauzet2016nitsche,hansbo2015characteristic,burman2014unfitted,mayer20103d,massing2015nitsche,gerstenberger2008extended,zonca2018unfitted}.
All of these methods sharply track the valve surfaces, but they entail possible issues regarding large mesh deformations and topological changes at valve closure \cite{spuhler20183d}.
On the other hand, Eulerian approaches, such as the immersed boundary method \cite{origIB,borazjaniCurvibFsiAv,griffithIb,bazilevsHughesIb,hsuImmersogeom,wangIbFem,zhangIbFem,nestola2019immersed}, the fictitious domain method \cite{glowinskiFd,gerbeauFd,bazilevsFd,bazilevsHughesFd,dehartFdAV,fdFsiRigidcontact,morsiAleTrileaflet} or the Resistive Immersed Implicit Surface (RIIS) method \cite{fedele2017patient,fumagalli2020image}, hinge upon an implicit representation of the leaflets and do not require mesh conformity between the fluid domain and the valves.
For further details and comparisons among different valve models, we refer the reader to \cite{reviewMarom,reviewVotta,hirschhorn2020fluid}.

In most of the abovementioned  {simulations}, however, the isovolumetric phases of the heartbeat are neglected due to the non-unique definition of pressure in the ventricle when all valves are closed \cite{schenkel2009mri,bavo2017patient,this2020augmented,quarteroni2010numerical,zingaro2022geometric}. This shortcoming is related to the absence of a stress condition on the fluid domain, that would otherwise ensure a correct description of the pressure during isovolumetric phases. This is observed for instance in \cite{bucelli2022partitioned,nordsletten2011fluid,viola2020fluid,zhang2001analysis, bucelli2022mathematical}, where fully coupled FSI models are used. 

Some studies have circumvented this issue by introducing a slight compressibility of blood -- see, e.g., \cite{doyle2015application,zhang2001analysis}.
However, this assumption may affect the simulation results also in the ejection and filling phases, and the assumption of blood incompressibility is quite established in the cardiovascular modeling community \cite{quarteroni2010numerical,incompressible1,incompressible2,formaggia2010cardiovascular}.
A way to overcome pressure non-determination, while preserving incompressibility, is provided by the Augmented Resistive Immersed Surface (ARIS) proposed in \cite{this2020augmented}: when both the mitral and the aortic valves are closed, the RIS method  {is augmented with a source term} concentrated on the valves, to impose a prescribed value for the pressure.

In this work, we  {introduce} an \emph{Augmented Resistive Immersed Implicit Surface} (ARIIS) method that extends the ARIS capability of treating isovolumetric phases to the framework of the RIIS method, thus supporting a mesh that is not conforming with the valves (cf.~\cref{tab:ris-riis-aris-ariis}). 
To quantitatively assess the results of the method, we propose a simulation setting in a simplified geometry that retains all the characteristics of the heart cycle and may be employed as a benchmark for cardiac hemodynamic solvers.
Moreover, we discuss the application of our method to a realistic geometry of the left heart, with a prescribed displacement coming from electromechanical simulations.

The structure of the paper is the following.
In \cref{sec:model}, we recall the RIIS method and derive the ARIIS method to prescribe the intraventricular pressure.
Then, in \cref{sec:results}, we assess our new method
in different scenarios: first, in \cref{sec:toyproblem}, we analyze the idealized case discussed in \cite{this2020augmented}; then, in \cref{sec:cardiocylinder}, we propose a simplified benchmark setting entailing ventricular contraction; finally, a cardiac case in a realistic geometry is considered in \cref{sec:cfd-lh}.

\begin{table}
\centering
\begin{tabular}{c|cc}
% \hline
& conforming mesh & non-conforming mesh
\\
\hline
no isovolumetric phases &RIS \cite{astorino2012robust} & RIIS \cite{fedele2017patient} \\
isovolumetric phases & ARIS \cite{this2020augmented}& ARIIS   
\\
% \hline
\end{tabular}
\caption{Features characterizing the RIS, RIIS, ARIS and ARIIS methods.}
\label{tab:ris-riis-aris-ariis}
\end{table}

\section{Mathematical model}\label{sec:model}
In this section, we describe  {the cardiac hemodynamic model} and we introduce a new augmented version of the RIIS method. Specifically, \Cref{sec:RIIS} is devoted to the Navier-Stokes equations in ALE framework with RIIS modeling of valves, and \Cref{sec:ARIIS} to the derivation of the ARIIS method.

%%%%%%%%%%%%%%%%%%%%%%%%%%%%%%%%%%%%%%%%%%%%%%%%%%%%%%%%%%%%%%%%%%%%%%
\subsection{The RIIS method for Navier-Stokes equations in ALE form}\label{sec:RIIS}
%%%%%%%%%%%%%%%%%%%%%%%%%%%%%%%%%%%%%%%%%%%%%%%%%%%%%%%%%%%%%%%%%%%%%%

In heart chambers, blood can be considered as an incompressible, viscous and Newtonian fluid \cite{quarteroni2019mathematical}. Let $\mathbf u:\Omega_t \times (0, T) \to \mathbb R^3$ and $p:\Omega_t \times (0, T) \to \mathbb R$ be the fluid velocity and pressure, respectively, where $T$ is the final computational time, and $\Omega_t$ the domain in current configuration at time $t$, with $t \in (0, T)$. The domain at any time $t$ is defined in terms of a displacement field $\mathbf d:\Omega_0 \times (0, T) \to \mathbb R^3$ as follows:
\[
    \Omega_t = \left\{ \mathbf x \in \mathbb{R}^3 : \mathbf x = \mathbf x_0 + \mathbf d(\mathbf x_0, t), \mathbf x_0 \in \Omega_0 \right\}\;.
\]
Furthermore, we denote by $\mathbf u_\ale :\Omega_t \times (0, T) \to \mathbb R^3$ the ALE velocity \cite{donea1982arbitrary,hughes1981lagrangian} and we compute it by deriving $\mathbf d$ with respect to time. The domain displacement is the solution of the following harmonic extension problem: 
\begin{subnumcases}{\label{eq:lifting}}
	- \div (K \grad \mathbf d) = \mathbf 0 &  in $\Omega_0 \times (0, T),$ \\
	\mathbf d = \mathbf d_{\partial \Omega}(\mathbf x, t) &  on $\partial \Omega_0 \times (0, T),$
\end{subnumcases}
where $\mathbf d_{\partial \Omega}:\partial \Omega_0 \times (0, T) \to \mathbb R^3$ is the boundary displacement (which is prescribed), and $K$ is a stiffness tensor field introduced to avoid distortion of mesh elements. 

To model the cardiac valves with the RIIS method, we consider a surface $\Gamma_\k$ immersed in $\Omega_t$, with $\k \in \mathcal I_{\mathrm{v}}$ (the set of immersed surfaces). We impose kinematic coupling between the surface and the fluid by penalizing the mismatch between the relative fluid velocity $\mathbf u - \mathbf u_\ale$ and the velocity of the immersed surface $ \mathbf u_{\Gamma_\k}$. 
Each surface is implicitly described by a signed distance function  $\varphi_\k: \Omega_t \times (0, T) \to \mathbb R$, such that $\Gamma_\k  = \left \{ \mathbf{x} \in \Omega_t \, : \, \varphi_\k ({\mathbf{x}}) = 0 \right \} $, for all $\k \in \mathcal I_\mathrm{v}$.  $\Gamma_\k$ is characterized by a resistance coefficient $R_\k$ and a parameter $\varepsilon_\k$ representing the half-thickness of the valve. The penalization is imposed in a narrow layer around $ \Gamma_\k $, represented by the following smoothed Dirac delta function:
 \begin{equation*} 
\delta_\k (\varphi_\k({\mathbf{x}})) = 
\begin{cases} 
	\dfrac{1 + \cos(\pi \varphi_\k({\mathbf{x}}) / \varepsilon_\k)}{2 \varepsilon_\k} & \text{ if } \;|\varphi_\k({\mathbf{x}})| \leq \varepsilon_\k, \\ 0 & \text{ if } \;|\varphi_\k({\mathbf{x}})| > \varepsilon_\k, \\ 
\end{cases} 
\end{equation*} 
  {with $\mathbf x \in \Omega_t$ and for all $\k \in \mathcal I_\mathrm{v}$}.
For additional details on the RIIS method, we refer the reader to \cite{fedele2017patient, fumagalli2020image}. 

The incompressible Navier-Stokes equations in the ALE framework with RIIS modeling of cardiac valves read as follows \cite{fumagalli2020image}:
\begin{subnumcases}{}
	\rho\left(\pdv{\mathbf u}{t} + \left(\left(\mathbf u - \mathbf u_\ale\right)\cdot\grad\right) \mathbf u\right) - \div\left(\mu (\grad \mathbf u + \grad^T \mathbf u ) \right)  + \grad p
	+ \sum_{\mathrm k \in \mathcal I_{\mathrm v}}\frac{R_\mathrm{k}}{\varepsilon_\mathrm{k}}\delta_\k(\varphi_\mathrm{k})(\mathbf u - \mathbf u_\ale - \mathbf u_{\Gamma_\k}) = \mathbf 0
	& in $\Omega_t \times (0, T)$, \label{eq:ns-momentum} 
	\\
	\div\mathbf{u} = 0 & in $\Omega _t \times (0, T)$, \label{eq:ns-continuity}
\end{subnumcases}
endowed with suitable initial and boundary conditions. We denote the different terms appearing in \eqref{eq:ns-momentum} as follows:
\begin{itemize}
	\item inertial term: $\bm {\bm{\mathcal I}}(\mathbf u) = \rho\left(\pdv{\mathbf u}{t} + \left(\left(\mathbf u - \mathbf u_\ale\right)\cdot\grad\right) \mathbf u\right)$;
	\item viscous term: $\bm{\mathcal D}(\mathbf u) = \div\left(\mu (\grad \mathbf u + \grad^T \mathbf u ) \right)$;
	\item resistive term: $\bm{\mathcal R}(\mathbf u) = \sum_{\mathrm k \in \mathcal I_{\mathrm v}}\frac{R_\mathrm{k}}{\varepsilon_\mathrm{k}}\delta_\k(\varphi_\mathrm{k})(\mathbf u - \mathbf u_\ale - \mathbf u_{\Gamma_\k})$.
\end{itemize}

%%%%%%%%%%%%%%%%%%%%%%%%%%%%%%%%%%%%%%%%%%%%%%%%%%%%%%%%%%%%%%%%%%%%%%
\subsection{The ARIIS method}\label{sec:ARIIS}
In this section, we derive the ARIIS method starting from the equations of the fluid model. 
To keep the notation light, we drop henceforth the explicit dependence on time of the domain and its subsets.  

\begin{figure}
	\centering
	\includegraphics[width=0.6\textwidth]{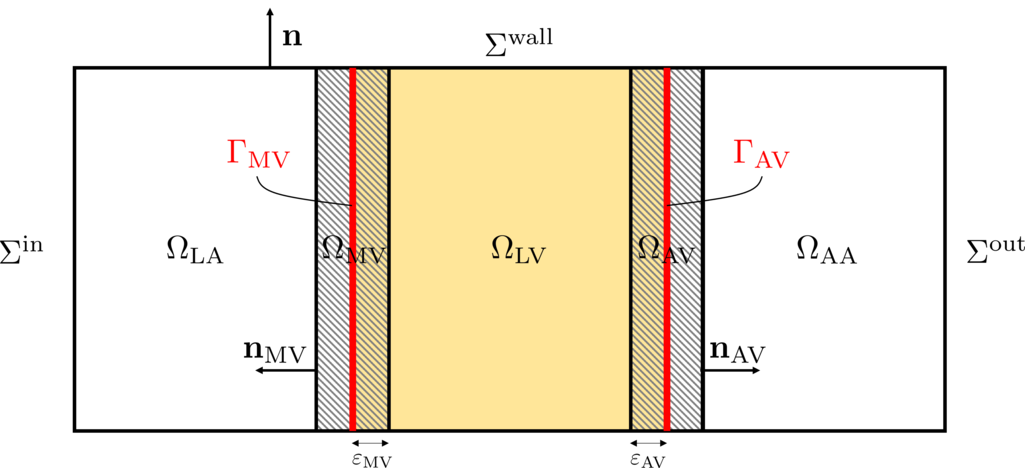}
	\caption{Sketch of the three-chambers domain $\Omega$ with its  {subsets} and boundaries. In yellow $\Omega_\lv$, in white $\Omega_\la,\Omega_\ao$, in striped pattern (partially overlapping $\Omega_\lv,\Omega_\la,\Omega_\ao$) the valve regions $\Omega_\mv,\Omega_\av$. The latter are defined by the immersed surfaces $\Gamma_\mv$ and $\Gamma_\av$ (in red) and the half-thicknesses $\varepsilon_\k$, with $\k \in \{ \mv, \av\}$. }
	\label{fig:domain_derivation}
\end{figure}

The left heart can be schematically outlined as a three-chambers domain as sketched in \Cref{fig:domain_derivation}:   the Left Atrium (LA) $\Omega_\la$,  the Left Ventricle (LV) $\Omega_\lv$ and the Ascending Aorta (AA) $\Omega_\ao$. These chambers are separated by two surfaces representing the Mitral Valve (MV) $\Gamma_\mv$ and the Aortic Valve (AV) $\Gamma_\av$, thus $\mathcal I_{\mathrm v} = \{ \mv, \av \}$. 

We denote by $\Omega$ the whole domain, such that $\overline{\Omega} = \overline{\Omega_\la} \cup \overline{\Omega_\lv} \cup \overline{\Omega_\ao}$. The domain boundary $\partial\Omega$ is partitioned into the inlet section $\Sigma^\mathrm{in}$, the outlet section $\Sigma^\text{out}$ and the wall $\Sigma^\mathrm{wall}$, as shown in \Cref{fig:domain_derivation}.
We introduce the sets
\begin{equation}
	\Omega_\k = \left\{\mathbf x \in \Omega \colon \mathrm{dist}(\mathbf x, \Gamma_\k) = \min_{\mathbf y\in\Gamma_\k}\|\mathbf x -\mathbf y\|<\varepsilon_\k \},\qquad \k\in\{\mv,\av\right\},
\end{equation}
where $\varepsilon_\k$ is the half thickness of $\Omega_\k$ (characterizing the RIIS method and already introduced in \Cref{sec:RIIS}), banded in \Cref{fig:domain_derivation}, for $\k\in\{\mv,\av\}$.
These regions have nontrivial intersections with the chambers defined above.

\begin{figure}
	\centering
	\includegraphics[width=0.4\textwidth]{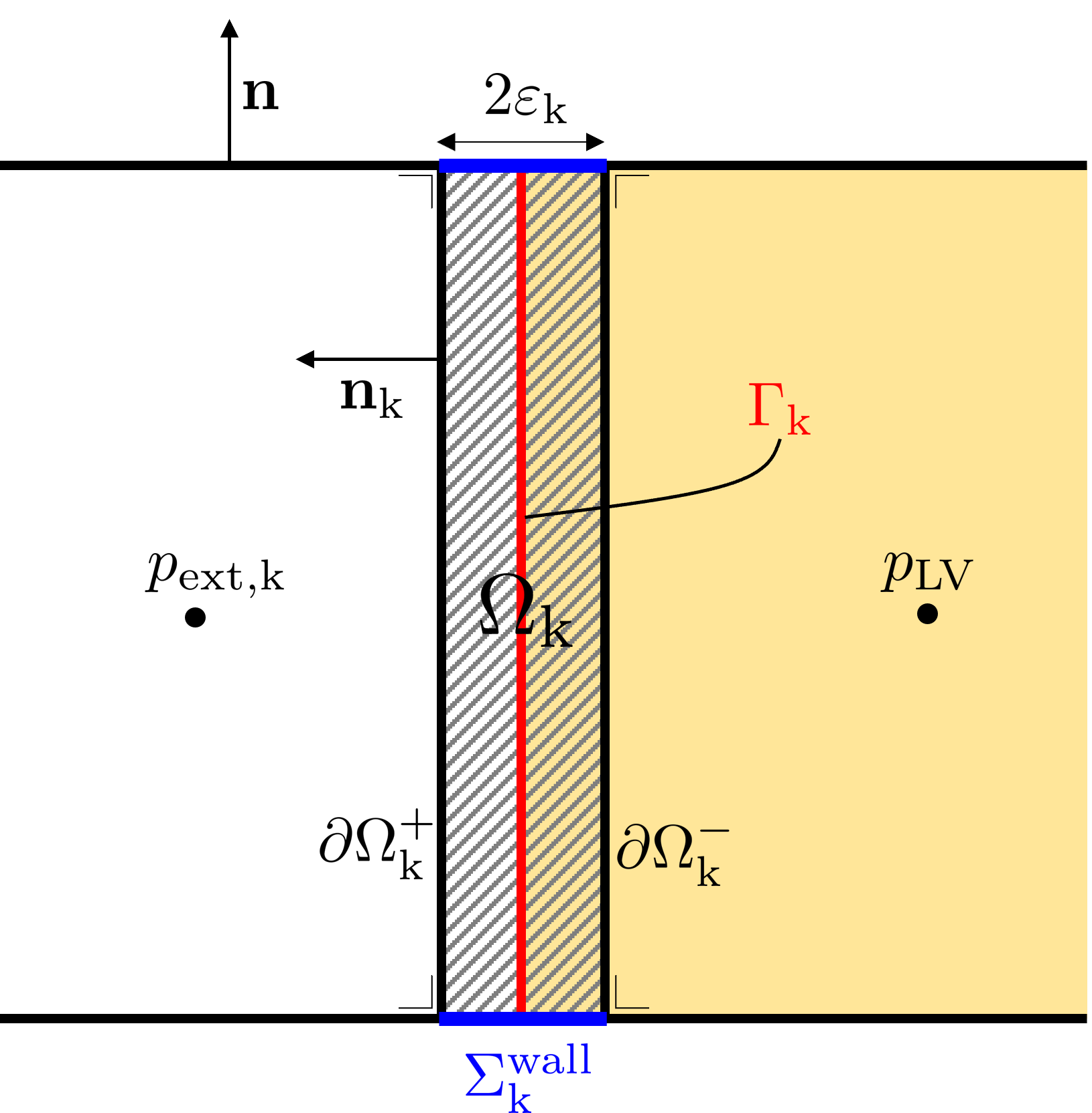}
	\caption{Sketch of the immersed surface $\Gamma_\k$ (red) with the corresponding valve region $\Omega_\k$ (striped pattern) and its boundaries.}
	\label{fig:valve-straight}
\end{figure}

With reference to \Cref{fig:valve-straight}, we introduce two geometric assumptions that will be used in the derivation of the augmented method of \Cref{sec:ARIIS}.

\begin{hp}(Flat valve surfaces)\label{hp:zero-curvature}
	For $\k\in\{\mv,\av\}$, the normal vector $\mathbf n_\k$ to the valve surface $\Gamma_\k$ (pointing outwards w.r.t.~$\Omega_\lv$) is constant over $\Gamma_\k$.
\end{hp}

\begin{rem} \label{rem:geometric}
 {If} Assumption \ref{hp:zero-curvature} is satisfied,  {we can define a constant vector field extending the definition of the valve normal vector $\mathbf n_\k$ to the whole valve region $\Omega_\k$ }. We denote such field with the same symbol $\mathbf n_\k:\Omega_\k\to\mathbb R^3,\k\in\{\mv,\av\}$.
\end{rem}

\begin{hp}(Valves orthogonal to the wall)\label{hp:normal-wall}
	 {By denoting with} $\mathbf n$ the normal vector of $\partial \Omega$, $\mathbf n_\k\cdot\mathbf n = 0$ on $\Sigma^\mathrm{wall}_\k$, for $\k\in\{\mv,\av\}$, where $\Sigma^\mathrm{wall}_\k=\Sigma^\mathrm{wall}\cap\partial\Omega_\k$.
\end{hp}

\begin{rem} \label{rem:areas}
 {By introducing} $\partial\Omega_\k^-=\partial\Omega_\k\cap\Omega_\lv$ and $\partial\Omega_\k^+=\partial\Omega_\k\setminus(\partial \Omega_\k^- \cup\Sigma^\mathrm{wall}_\k)$ (cfr.~\Cref{fig:valve-straight}) we observe that $|\partial\Omega_\k^-|=|\partial\Omega_\k^+|=|\Gamma_\k|,\,\k\in\{\mv,\av\}$.
\end{rem}

Moreover, we make the following assumptions:
\begin{hp}[Constant pressure in the compartments]
	Pressure is constant in space within $\Omega_\mathrm{LA}$, $\Omega_\mathrm{LV}$ and $\Omega_\ao$. We will denote the respective constant values with $p_\mathrm{LA}(t)$, $p_\mathrm{LV}(t)$ and $p_\ao(t)$.
	\label{hp:constant-pressure}
\end{hp}
\begin{hp}[Negligible inertia and viscosity within valves]
	For $\mathrm{k} \in \{\text{MV}, \text{AV}\}$, inertial and viscous terms in \eqref{eq:ns-momentum} are negligible in $\Omega_\k$: $\bm{\mathcal I}(\mathbf u) \approx \mathbf 0$ and $\bm{\mathcal D}(\mathbf u) \approx \mathbf 0$.
	\label{hp:negligible-inertia-viscosity}
\end{hp}

 {When MV and AV are closed, the intraventricular pressure is prone to spurious oscillations, due to the ventricle being fully enclosed by boundaries on which a Dirichlet-type condition on the velocity is imposed (either strongly or through the RIIS penalty term). Thus, we augment \Cref{eq:ns-momentum} with a reaction term to impose $p^*: (0, T)\to \mathbb R$, a prescribed value of the ventricular pressure (constant in space by Assumption~\ref{hp:constant-pressure}), that can be obtained, for instance, from an electromechanical simulation or from patient-specific measured data.}

We assume the perturbation term to be in the form:
\begin{equation}
	 \sum_{\k \in \{\mv, \av \}} \mathcal C_\k \delta_\k \mathbf n_\k,
	 \label{eq:perturbation}
\end{equation}
with $\mathcal C_\k \in \mathbb{R}$, for $\k \in \{ \mv, \av \}$. 
This choice of the reaction term \eqref{eq:perturbation} is such that the augmented formulation acts on the valves only and does not perturb the momentum equation outside $\Omega_\k$. 

The pertubation term represents the force that the blood exerts on the closed valves during isovolumetric phases. We derive it to enforce that the ventricular pressure matches the reference one $p^*$. Thus, following \cite{this2020augmented}, we derive an estimation of the ventricular pressure $p_\lv(t)$ when both valves are closed. The estimate will be used to determine the corrective term $\mathcal C_\k$ in \eqref{eq:perturbation}.

From \eqref{eq:ns-momentum} and \Cref{hp:negligible-inertia-viscosity}, we deduce
\[
\grad p + \bm{\mathcal R}(\mathbf u) = \mathbf 0 \qquad \text{in } \Omega_\k,
\]
for all $\mathrm{k} \in \{\text{MV}, \text{AV}\}$. Multiplying by $\mathbf n_\k$ and integrating over $\Omega_\k$, we get
\[
\int_{\Omega_\k}\left(\grad p + \bm{\mathcal R}(\mathbf u)\right)\cdot\mathbf n_\k = 0.
\]
By \Cref{hp:zero-curvature}, we can take $\mathbf n_\k$ out of the integral and integrate by parts the pressure term yielding
\begin{equation*}
	\left(\int_{\partial\Omega_\k}p\mathbf n + \int_{\Omega_\k}\bm{\mathcal R}(\mathbf u)\right)\cdot\mathbf n_\k = 0, 
	\end{equation*}
	\begin{equation*}
	\int_{\partial\Omega_\k}p\mathbf n\cdot\mathbf n_\k + \int_{\Omega_\k}\bm{\mathcal R}(\mathbf u)\cdot \mathbf n_\k = 0.
\end{equation*}
Using Assumptions \ref{hp:zero-curvature} and \ref{hp:normal-wall}, we get 
\begin{equation}
	\left(p_\mathrm{LV} - p_\mathrm{ext, k}\right) |\Gamma_\k| - \int_{\Omega_\k}\bm{\mathcal R}(\mathbf u)\cdot\mathbf n_\k = 0,
	\label{eq:pressure-estimate-single-surface}
\end{equation}
where $p_\mathrm{ext, k} = p_\mathrm{LA}$ for $\mathrm k = \mathrm{MV}$ and $p_\mathrm{ext, k} = p_\ao$ for $\mathrm k = \mathrm{AV}$. Finally, summing \eqref{eq:pressure-estimate-single-surface} for both valves, we obtain:
\begin{equation}
	\left(p_\mathrm{LV} - p_\mathrm{LA}\right) |\Gamma_\mathrm{MV}| + \left(p_\mathrm{LV} - p_\ao\right) |\Gamma_\mathrm{AV}| -
	\sum_{\k \in \{ \mathrm{MV}, \mathrm{AV}\}}  \int_{\Omega_\k}\bm{\mathcal R}(\mathbf u)\cdot\mathbf n_\k = 0,
	\label{eq:pressure-estimate-both-surfaces-derivation}
\end{equation}
from which we derive
\begin{equation}
	p_\mathrm{LV} = \frac{1}{|\Gamma_\mathrm{MV}| + |\Gamma_\mathrm{AV}|}\left(|\Gamma_\mathrm{MV}|p_\mathrm{LA} + |\Gamma_\mathrm{AV}|p_\ao + \sum_{\k \in \{ \mathrm{MV}, \mathrm{AV}\}}  \int_{\Omega_\k}\bm{\mathcal R}(\mathbf u)\cdot\mathbf n_\k\right)\;.
	\label{eq:pressure-estimate-both-surfaces}
\end{equation}

Repeating these calculations including the perturbation term \eqref{eq:perturbation}, \eqref{eq:pressure-estimate-both-surfaces-derivation} rewrites as
\[
\sum_{\k \in \{\mv, \av\}}\left ( \left(p_\mathrm{LV} - p_\mathrm{ext, k}\right)|\Gamma_\k| - \int_{\Omega_\k}\bm{\mathcal R}(\mathbf u)\cdot\mathbf n_\k - \mathcal{C}_\k\int_{\Omega_\k}\delta_\k \right ) = 0\;,
\]
so that, if the perturbation satisfies
\begin{equation}
	\sum_{\k \in \{\mv, \av\}} \int_{\Omega_\k}\mathcal{C}_\k\delta_\k = \sum_{\k \in \{\mv, \av\}} \left ( (p^* - p_\mathrm{ext, k})|\Gamma_\k| - \int_{\Omega_\k}\bm{\mathcal R}(\mathbf u)\cdot\mathbf n_\k \right ),
	\label{eq:condition-ck}
\end{equation}
then our estimate for $p_\mathrm{LV}$ becomes $p_\mathrm{LV} = p^*$.

Observing that $\int_{\Omega_\k}\delta_\k = |\Gamma_\k|$, a possible definition of the corrective term satisfying \eqref{eq:condition-ck} is:

\begin{equation}
	\mathcal{C}_\k  {(\mathbf u, p)} = p^* - p_\mathrm{ext, k} - \frac{1}{|\Gamma_\mathrm{MV}| + |\Gamma_\mathrm{AV}|} \sum_{\k \in \{ \mathrm{MV}, \mathrm{AV}\}}  \int_{\Omega_\k}\bm{\mathcal R}(\mathbf u)\cdot\mathbf n_\k\;.
	\label{eq:ck2}
\end{equation}

Thus, the ARIIS method consists in solving the following problem:
\begin{subnumcases}{\label{eq:ns-ariis}}
	\begin{aligned}
		& \rho\left(\pdv{\mathbf u}{t} + \left(\left(\mathbf u - \mathbf u_\ale\right)\cdot\grad\right) \mathbf u\right) - \div\left(\mu(\grad \mathbf u + \grad^T \mathbf u) \right) + \grad p\\
		& + \sum_{\mathrm k \in \mathcal \{\mv, \av \}} \left ( \frac{R_\mathrm{k}}{\varepsilon_\mathrm{k}}\delta_\k(\varphi_\mathrm{k})(\mathbf u - \mathbf u_\ale - \mathbf u_{\Gamma_\k})  + \chi_\mathrm{iso}(t) \mathcal{C}_\k  {(\mathbf u, p)} \delta _\k \mathbf n_\k \right )= \mathbf 0
	\end{aligned}
	& in $\Omega_t \times (0, T)$, \label{eq:ns-momentum-ariis} \\
	\div\mathbf{u} = 0 & in $\Omega _t \times (0, T)$, \label{eq:ns-continuity-ariis}
\end{subnumcases}
endowed with suitable initial and boundary conditions. $\chi_\mathrm{iso}(t)$ is a characteristic function equal to 1 during the isovolumetric phases, 0 otherwise: we activate the ARIIS correction term only when both valves are simultaneously closed. $\chi_\mathrm{iso}(t)$ can be prescribed a priori or be determined by pressure jump conditions (to determine the opening and closing of valves) \cite{this2020augmented}.  

%%%%%%%%%%%%%%%%%%%%%%%%%%%%%%%%%%%%%%%%%%%%%%%%%%%%%%%%%%%%%%%%%%%%%%
%%%%%%%%%%%%%%%%%%%%%%%%%%%%%%%%%%%%%%%%%%%%%%%%%%%%%%%%%%%%%%%%%%%%%%
\section{Numerical results}\label{sec:results}
%%%%%%%%%%%%%%%%%%%%%%%%%%%%%%%%%%%%%%%%%%%%%%%%%%%%%%%%%%%%%%%%%%%%%%
%%%%%%%%%%%%%%%%%%%%%%%%%%%%%%%%%%%%%%%%%%%%%%%%%%%%%%%%%%%%%%%%%%%%%%
In this section, we present and discuss the results on the ARIIS method by carrying out numerical simulations  {of} three different problems. All three tests feature valves that open and close. In \Cref{sec:toyproblem}, we check the validity of our method by considering the simple problem introduced in \cite{this2020augmented} (Test A). In \Cref{sec:cardiocylinder}, we propose a new benchmark problem consisting of the flow in a compliant pipe with ventricle-like shortening (Test B). Finally, in \Cref{sec:cfd-lh}, we apply our method to a cardiac case, i.e. the flow in a realistic left heart geometry (Test C).  

The physical parameters for  blood are density {$\rho = \SI{1.06e3}{\kilo\gram\per\metre\cubed} $} and dynamic viscosity $\mu=\SI{3.5e-3}{\kilo\gram\per\metre\per\second}$. In all the numerical experiments considered, we apply a null velocity initial condition. Furthermore, similarly to \cite{fumagalli2020image, zingaro2022geometric}, we use a quasi-static approach by choosing $\mathbf u_{\Gamma_\k} = \mathbf 0$. 

We discretize \eqref{eq:ns-ariis} in space with piecewise linear Finite Elements (FE) for velocity and pressure ($\mathbb P_1 - \mathbb P_1$) and in time with the backward Euler method. We employ a semi-implicit treatment of the non-linear term. In \Cref{sec:toyproblem} and \Cref{sec:cardiocylinder}, we use a
SUPG-PSPG stabilization \cite{TEZDUYAR19911}. Differently, in \Cref{sec:cfd-lh}, we use the
VMS-LES method acting as both a stabilization method and a turbulence model to account for the transition-to-turbulence flow regime typically occurring in cardiac flows \cite{bazilevs2007variational, forti2015semi, zingaro2021hemodynamics}. The lifting problem \eqref{eq:lifting} is discretized with linear FEs. 

We carry out our numerical simulations in \texttt{life$^\texttt{x}$} \cite{africa2022lifex}\footnote{\url{https://lifex.gitlab.io/}}, a high-performance \texttt{C++} FE library developed within the iHEART
project\footnote{iHEART - An Integrated Heart model for the simulation of the cardiac function, European Research Council (ERC) grant agreement No 740132, P.I. Prof. A. Quarteroni, 2017-2022}, mainly focused on cardiac simulations and based on the \texttt{deal.II} finite element core \cite{arndt2021dealii,arndt2020dealii,dealii}.

\subsection{Test A: a simple benchmark problem}
\label{sec:toyproblem}

In this section, following \cite{this2020augmented}, we consider a benchmark problem that was originally introduced to test the ARIS method in a simplified setting.

% Domain (spatial and temporal).
The domain is a cylinder of radius $R_\text{c} = \SI{0.01}{\metre}$ and length $L_\text{c} = \SI{0.1}{\metre}$. It is divided into three cylindrical compartments, representing, in an idealized context, the LA, LV and AA, of lengths $L_\la = \SI{0.02}{\metre}$, $L_\lv = \SI{0.06}{\metre}$ and $L_\ao = \SI{0.02}{\metre}$, respectively. Two planar surfaces represent the MV and AV. We solve  {in} the time interval $[0, T]$, with $T = \SI{0.2}{\second}$.

\begin{figure}
    \centering
    \begin{subfigure}{0.49\textwidth}
        \includegraphics[width=\textwidth]{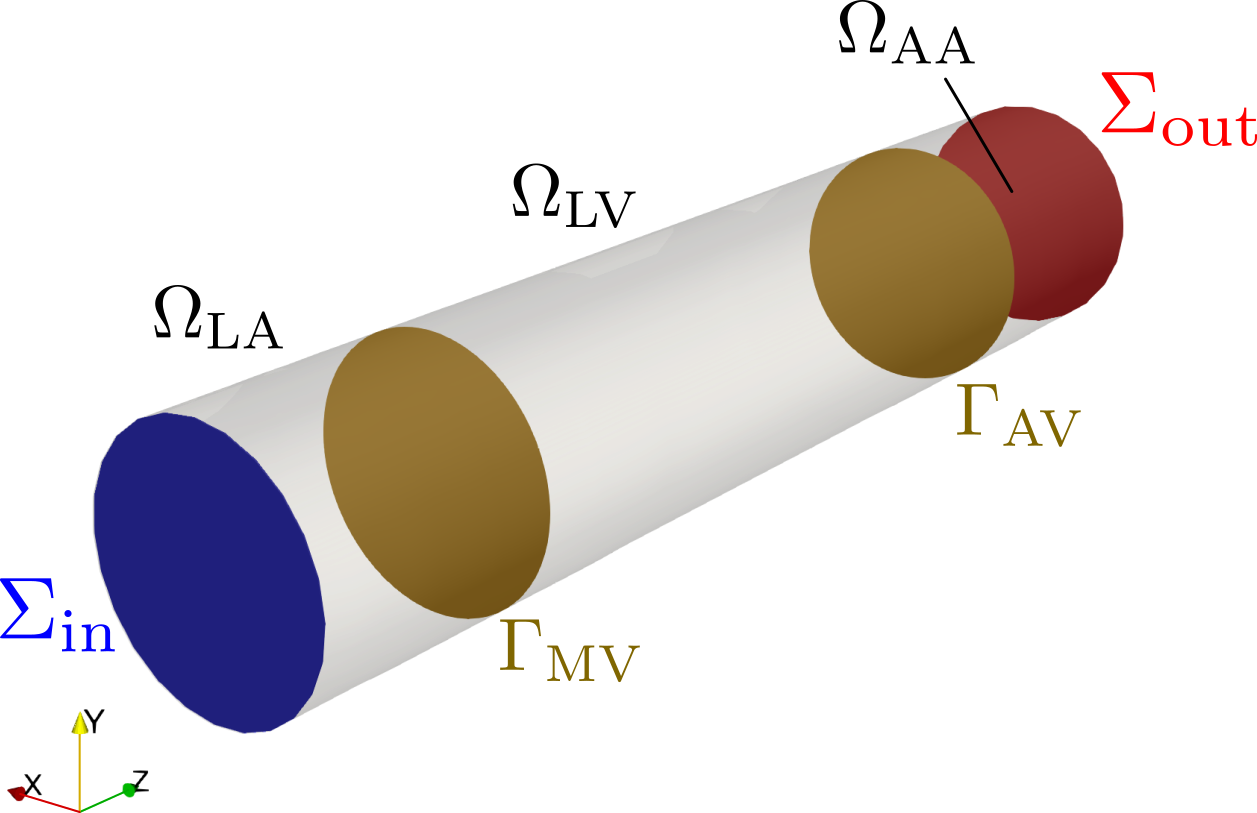}
        \caption{}
    \end{subfigure}
    \begin{subfigure}{0.49\textwidth}
        \centering
        \includegraphics[width=\textwidth]{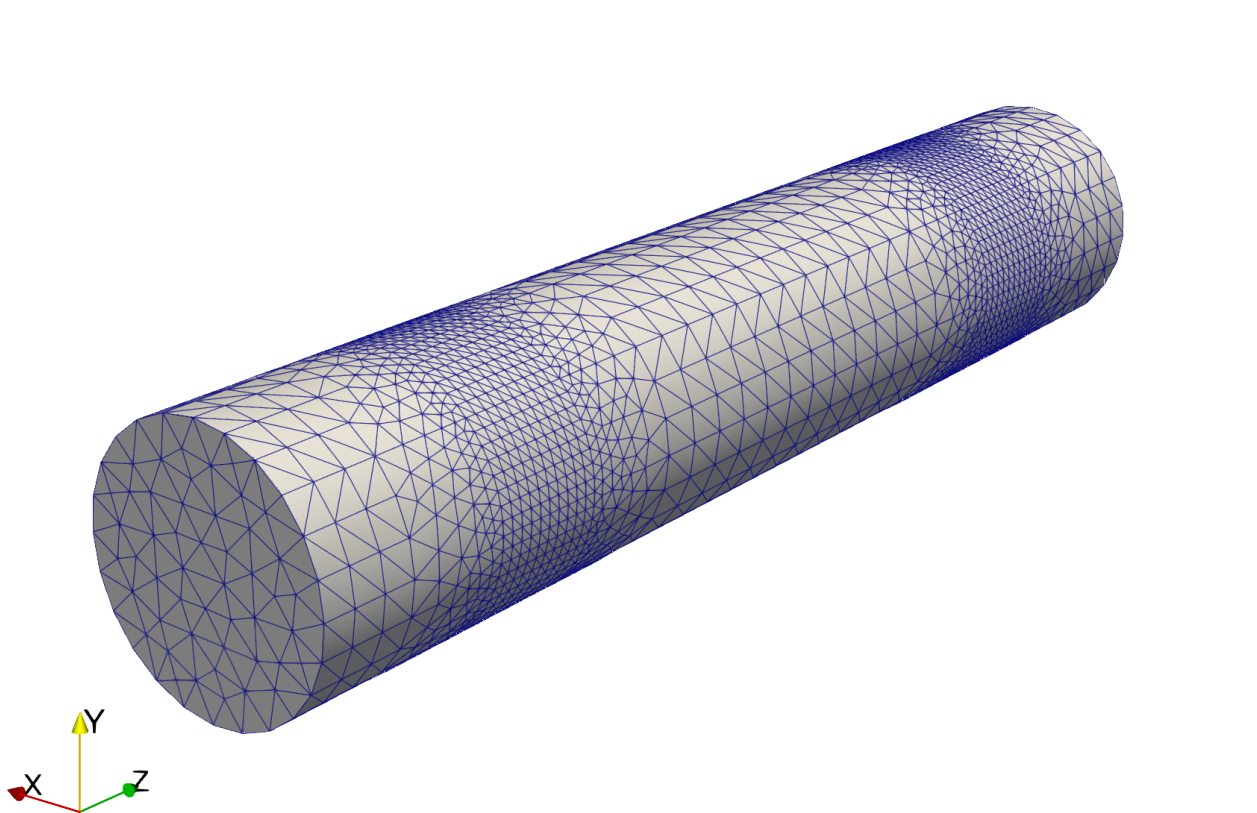}
        \caption{}
    \end{subfigure}
    
    \caption{Test A. Domain for the cylinder test cases, with highlighted immersed and boundary surfaces (a); tetrahedral mesh for the cylinder test cases.}
    \label{fig:cylinder-mesh}
\end{figure}

% Mesh.
The domain is discretized with a tetrahedral mesh of \num{75933} elements, for a total of \num{56684} degrees of freedom. The mesh is finer near to the immersed surfaces, to better capture their presence, with a minimum element diameter $h_{\min} = \SI{1}{\milli\metre}$ and a maximum diameter $h_{\max} = \SI{4.6}{\milli\metre}$ (see \cref{fig:cylinder-mesh}).
Simulations ran in parallel using 4 cores of a local workstation, each with an Intel Core i5-9600K@3.70GHz processor.

\begin{figure}
    \centering
    \begin{subfigure}{0.4\textwidth}
        \includegraphics[width=\textwidth]{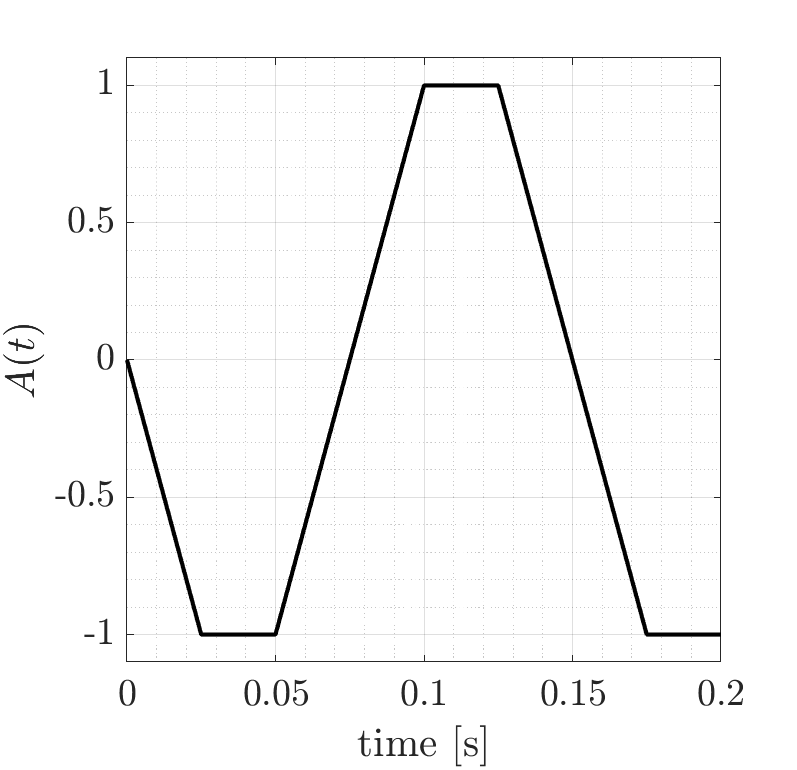}
        \caption{}
        \label{fig:this_displacement_factor}
    \end{subfigure}
    \begin{subfigure}{0.4\textwidth}
        \includegraphics[width=\textwidth]{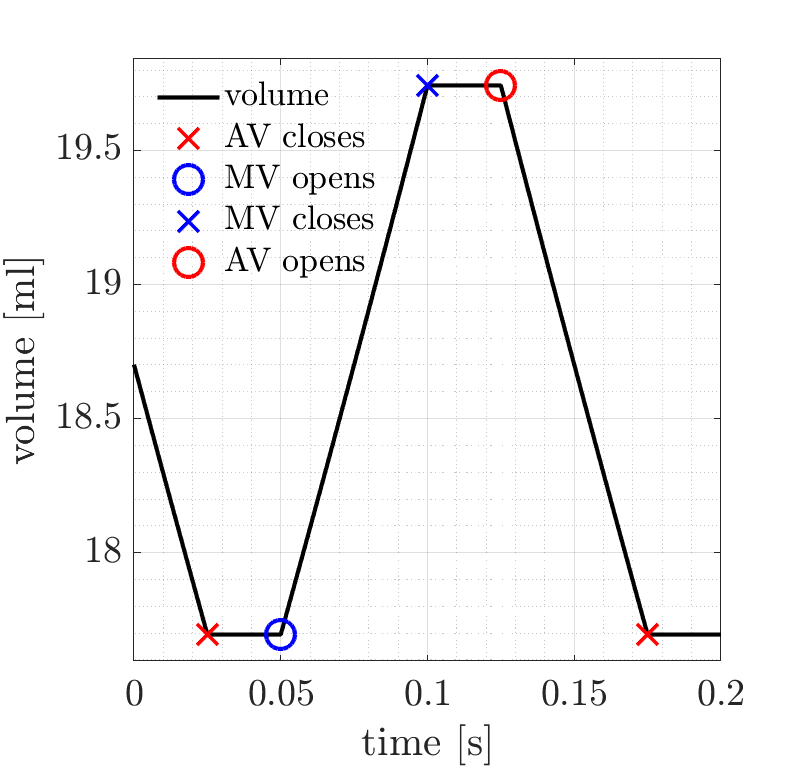}
        \caption{}
        \label{fig:this_volume}
    \end{subfigure}
    \caption{Test A. (a) Plot of the function $A(t)$ that defines the time evolution of the boundary displacement in the cylindrical toy problem. (b) The corresponding volume of the ventricular compartment with valve opening and closing times.}
\end{figure}

% Boundary conditions.
Following \cite{this2020augmented}, we impose a homogeneous and constant pressure of $p_\mathrm{in} = \SI{0}{\mmhg}$ at the inlet section, and a homogeneous and constant pressure of $p_\mathrm{out} = \SI{75}{\mmhg}$ at the outlet section. The displacement $\mathbf d_{\partial\Omega}$ of the lateral boundary is prescribed analytically and mimics the contraction-relaxation cycle of a human ventricle. For a given point $\mathbf x = (x_1, x_2, x_3)^T$ and time $t$, it is defined as
\begin{equation*}
    \mathbf d_{\partial\Omega}(\mathbf x, t) = \begin{cases}
        \overline{w} \, A(t) \,  \mathbf{e}_r(\mathbf x) \,  \exp\left(-\frac{\left|x_3 - \frac{L_\text{c}}{2}\right|^2}{2\sigma^2}\right) & \text{if } x_3 \in [L_\la, L_\la + L_\lv)\;, \\
        \mathbf 0 & \text{otherwise,}
    \end{cases}
\end{equation*}
with
\begin{equation*}
      {\mathbf{e}}_r(\mathbf x) = \frac{(x_1, x_2, 0)^T}{\sqrt{x_1^2 + x_2^2}}\;.
\end{equation*}
$A(t)$ is the piecewise linear function depicted in \Cref{fig:this_displacement_factor}. We set $\sigma = \SI{0.015}{\metre}$ and $\overline{w} = \SI{4.6e-4}{\metre}$, to have the same time evolution of volume as in \cite{this2020augmented} (see \cref{fig:this_volume}). We take $K = I$ in \eqref{eq:lifting}.

% Valve opening/closing.
We simulate the opening of a valve by instantaneously removing the corresponding surface from the domain. Valves are opened and closed at prescribed times, following the evolution of the volume of the ventricular compartment. The MV is closed when the simulation starts, while the AV is open. Closing and opening times are reported in \Cref{fig:this_volume}.

\begin{figure}
    \centering
    \includegraphics[width=0.7\textwidth]{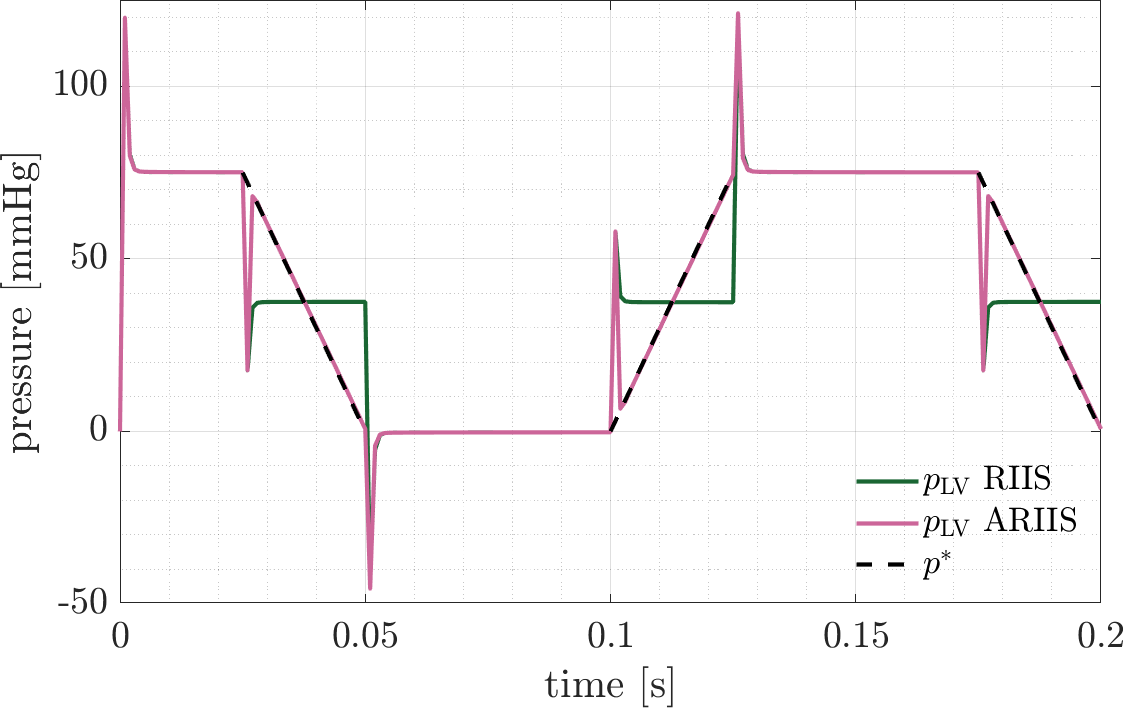}
    \caption{Test A. Ventricular and reference pressures for the cylindrical toy problem, with RIIS and ARIIS methods, using resistance $R = \SI{e4}{\kilo\gram\per\metre\per\second}$ and $\varepsilon = \SI{0.002}{\metre}$.}
    \label{fig:this_riis_vs_ariis}
\end{figure}

% Tests: RIIS vs ARIIS.
In this setting, we carry out a comparison of the results obtained with the RIIS method against those obtained with the ARIIS method,  {using as reference pressure $p^*(t)$ a piecewise linear function}. The evolution of ventricular pressure for both cases, computed with resistance $R = \SI{e4}{\kilo\gram\per\metre\per\second}$ and $\varepsilon = \SI{0.002}{\metre}$, is reported in \Cref{fig:this_riis_vs_ariis}. The plots show how the ARIIS method allows the ventricular pressure to accurately follow the provided reference pressure. The observed peaks are associated to the simplified and instantaneous way in which valves are opened and closed and to the explicit computation of the corrective term \eqref{eq:perturbation}.

\begin{figure}
    \centering
    
    \begin{subfigure}{0.594\textwidth}
        \includegraphics[width=\textwidth]{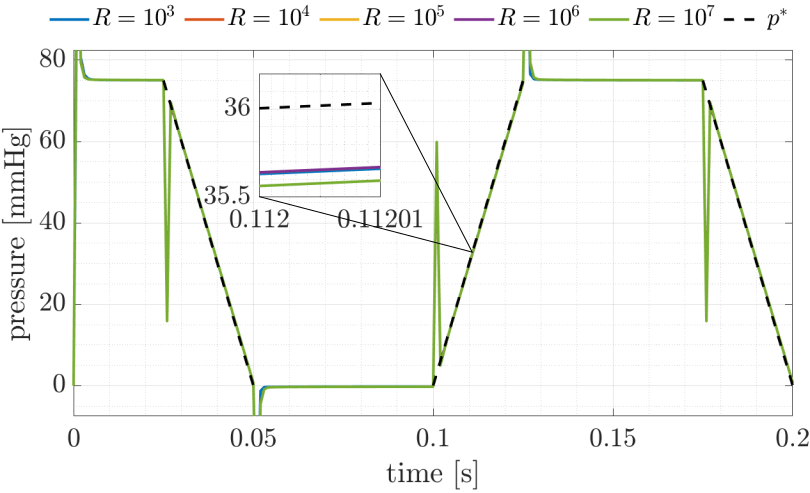}
        \caption{}
    \end{subfigure}
    \hspace{0.5em}
    \begin{subfigure}{0.27\textwidth}
        \includegraphics[width=\textwidth]{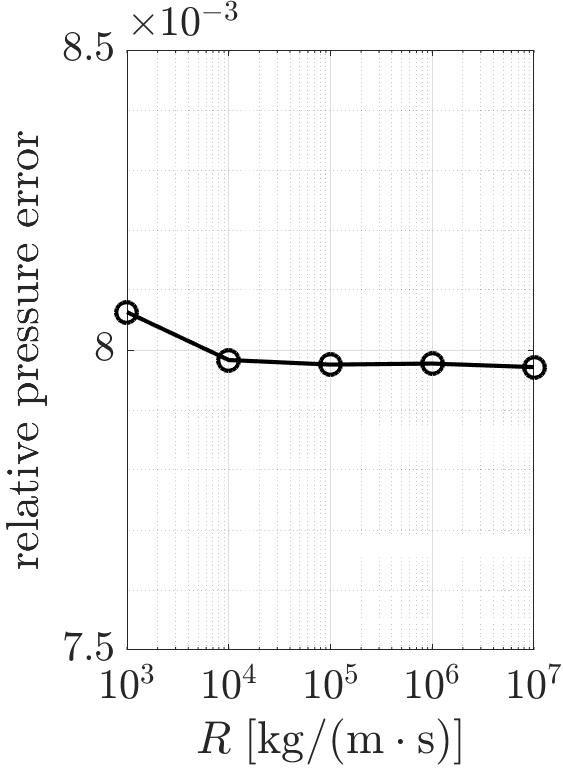}
        \caption{}
        \label{fig:this_pressure_error_vs_R}
    \end{subfigure}
    
    \caption{Test A. (a) Evolution of pressure with the ARIIS method with varying values of the resistance. (b) Relative error between the computed pressure during isovolumic phases and the reference pressure $p^*$.}
    \label{fig:this_resistance}
\end{figure}

% Tests: ARIIS with varying resistance.
Moreover, we carry out a sensitivity analysis by varying the resistance coefficient $R$ in the ARIIS method, to understand how the quality of the results is influenced by it. Results are reported in \Cref{fig:this_resistance}. Although the resistance coefficient varies by several orders of magnitude, no difference is observed on the accuracy of the ventricular pressure. This is evident in particular in \Cref{fig:this_pressure_error_vs_R}, reporting the relative pressure error
\begin{equation}
    \frac{\max_{t \in T_\text{iso}}|p_\lv - p^*|}{\max_{t \in T_\text{iso}}|p^*|}\;,
    \label{eq:relative-pressure-error}
\end{equation}
where $T_\text{iso} = \{ t \in (0, T) : \chi_\text{iso}(t) = 1\}$ is the set of times at which both valves are closed. The error is approximately equal to \num{8e-3} regardless of the value of $R$. The ARIIS method, therefore, yields reliable pressure results also with high values of $R$, that ensure negligible spurious flow through the resistive surfaces.

Overall, the obtained results indicate that the ARIIS method is successful in its aim of producing a ventricular pressure that closely follows the prescribed reference evolution.

\FloatBarrier

\subsection{Test B:  {a benchmark problem including ventricular shortening}}
\label{sec:cardiocylinder}

\begin{figure}
    \centering
    \includegraphics[width=0.6\textwidth]{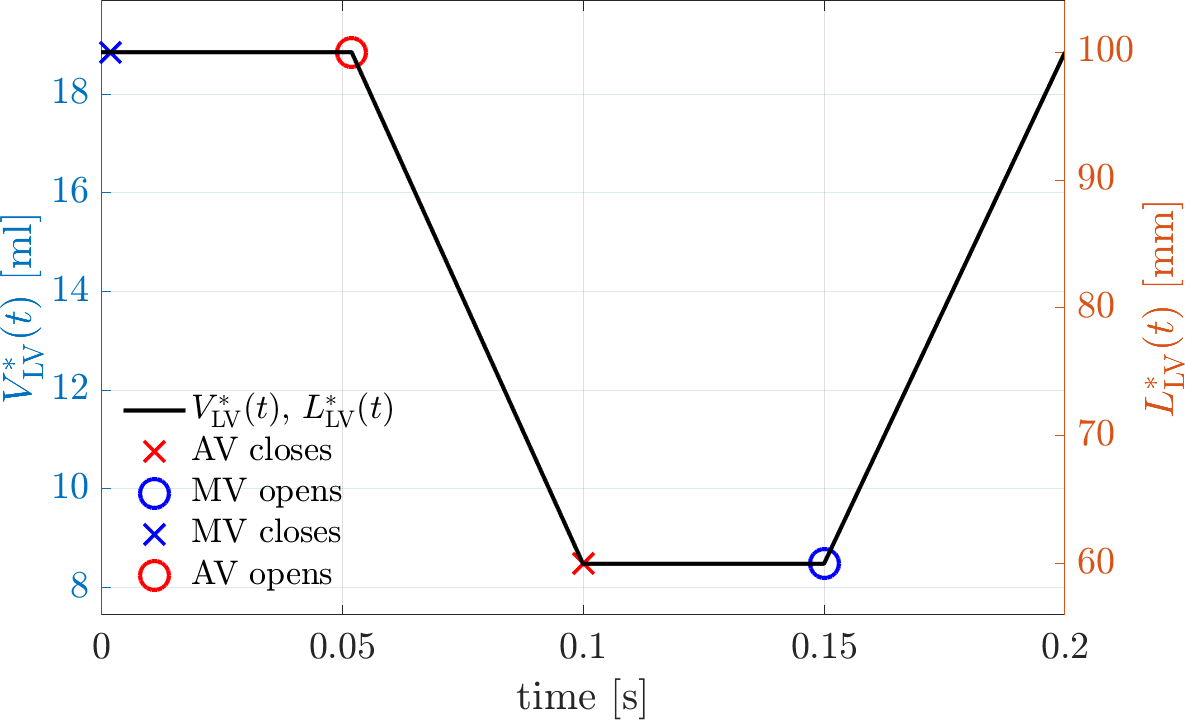}
    \caption{Test B. Prescribed ventricular volume $V_\lv^*(t)$ (left axis) and length $L_\lv^*(t)$ (right axis).}
    \label{fig:cardiocylinder-length-volume}
\end{figure}

% Changes in setup with respect to previous test.
As an intermediate step towards cardiac simulations, we introduce a novel test case in a cylindrical domain that mimics the ventricular shortening during contraction.  We use the same domain as in \Cref{sec:toyproblem}, but change the boundary displacement as follows:
\begin{align}
    \label{eq:displacement-cardiocylinder}
    \mathbf d_{\partial\Omega}(\mathbf x, t) &= \begin{cases}
        \mathbf 0 & \text{if } x_3 \in [0, L_\la)\;, \\ 
        \mathbf d_{\partial\Omega}^r(\mathbf x, t) + \mathbf d_{\partial\Omega}^z(\mathbf x, t) & \text{if } x_3 \in [L_\la, L_\la + L_\lv)\;, \\
        \mathbf (0,\; 0,\; L_\lv^*(t) - L_\lv)^T & \text{if } x_3 \in [L_\la + L_\lv, L)\;,
    \end{cases} \\
\end{align}
with
\begin{align}
    \mathbf d_{\partial\Omega}^r(\mathbf x, t) &= \left(R_\text{c} + c(t) \sin\left(\frac{\pi (x_3 - L_\la)}{L_\lv}\right)\right) \mathbf r(\mathbf x) - \mathbf x\;, \\
    \mathbf d_{\partial\Omega}^z(\mathbf x, t) &= \frac{x_3 - L_\la}{L_\lv}(L_\lv^*(t) - L_\lv)\;
    \end{align}
and
\begin{align}
    c(t) &= \frac{4R_\text{c}}{\pi} + \frac{\sqrt{16\,R_\text{c}^2\,L_\lv^*(t)^2 - 2\,\pi \, L_\lv^*(t)\, \left(\pi \, L_\lv^*(t) \, R_\text{c}^2 - V_\lv^*(t)\right)}}{\pi \, L_\lv^*(t)}\;.
\end{align}
In the above, $L_\lv^*(t)$ and $V_\lv^*(t)$ are prescribed  {time dependent functions} for the ventricular length and volume, respectively. The displacement is such that, at any time $t$, the ventricular length and volume in the deformed configuration match the prescribed ones. We take $K = I$ in \eqref{eq:lifting}. Valve positions change over time following the domain displacement. Their opening and closing times are reported in \Cref{fig:cardiocylinder-length-volume}. The MV starts open, and the AV starts closed. Moreover, we set inlet and outlet boundary conditions to $p_\text{in} = \SI{0}{\mmhg}$ and $p_\text{out} = \SI{80}{\mmhg}$, to replicate the typical range that characterizes the heart function.

Numerical simulations are run in parallel on the GALILEO100 supercomputer\footnote{528 computing nodes each 2 x CPU Intel CascadeLake 8260, with 24 cores each, 2.4 GHz, 384GB RAM. See \url{https://wiki.u-gov.it/confluence/display/SCAIUS/UG3.3\%3A+GALILEO100+UserGuide} for technical specifications.} at the CINECA supercomputing center, using 48 cores.

\begin{figure}
    \centering
    \begin{subfigure}{0.49\textwidth}
        \centering
        \textbf{RIIS}
    \end{subfigure}
    \begin{subfigure}{0.49\textwidth}
        \centering
        \textbf{ARIIS}
    \end{subfigure}
    
    \vspace{1em}
    \begin{subfigure}{\textwidth}
        \includegraphics[width=0.49\textwidth,
                         trim={0in 7in 0in 7in},
                         clip]{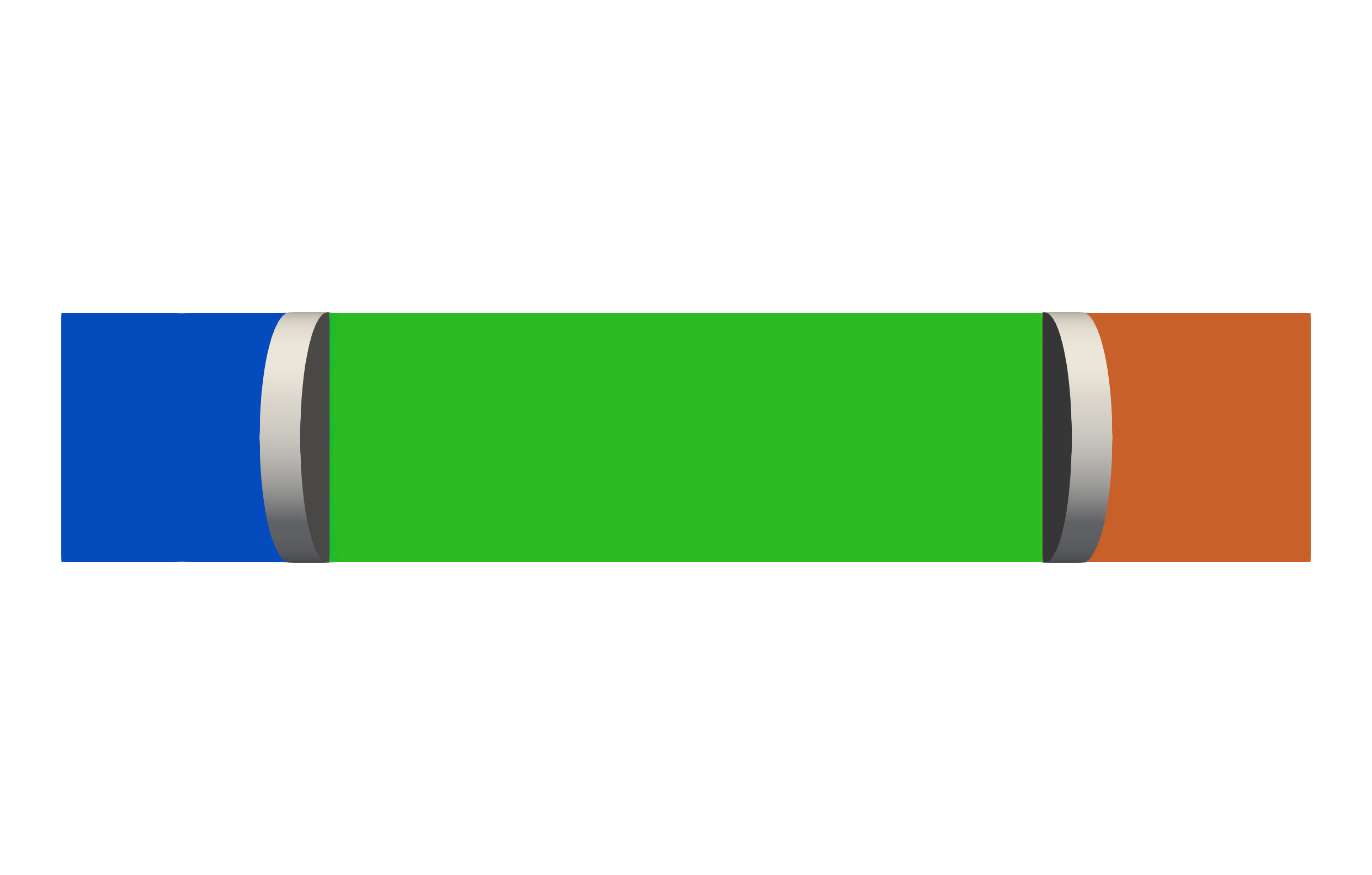}
        \includegraphics[width=0.49\textwidth,
                         trim={0in 7in 0in 7in},
                         clip]{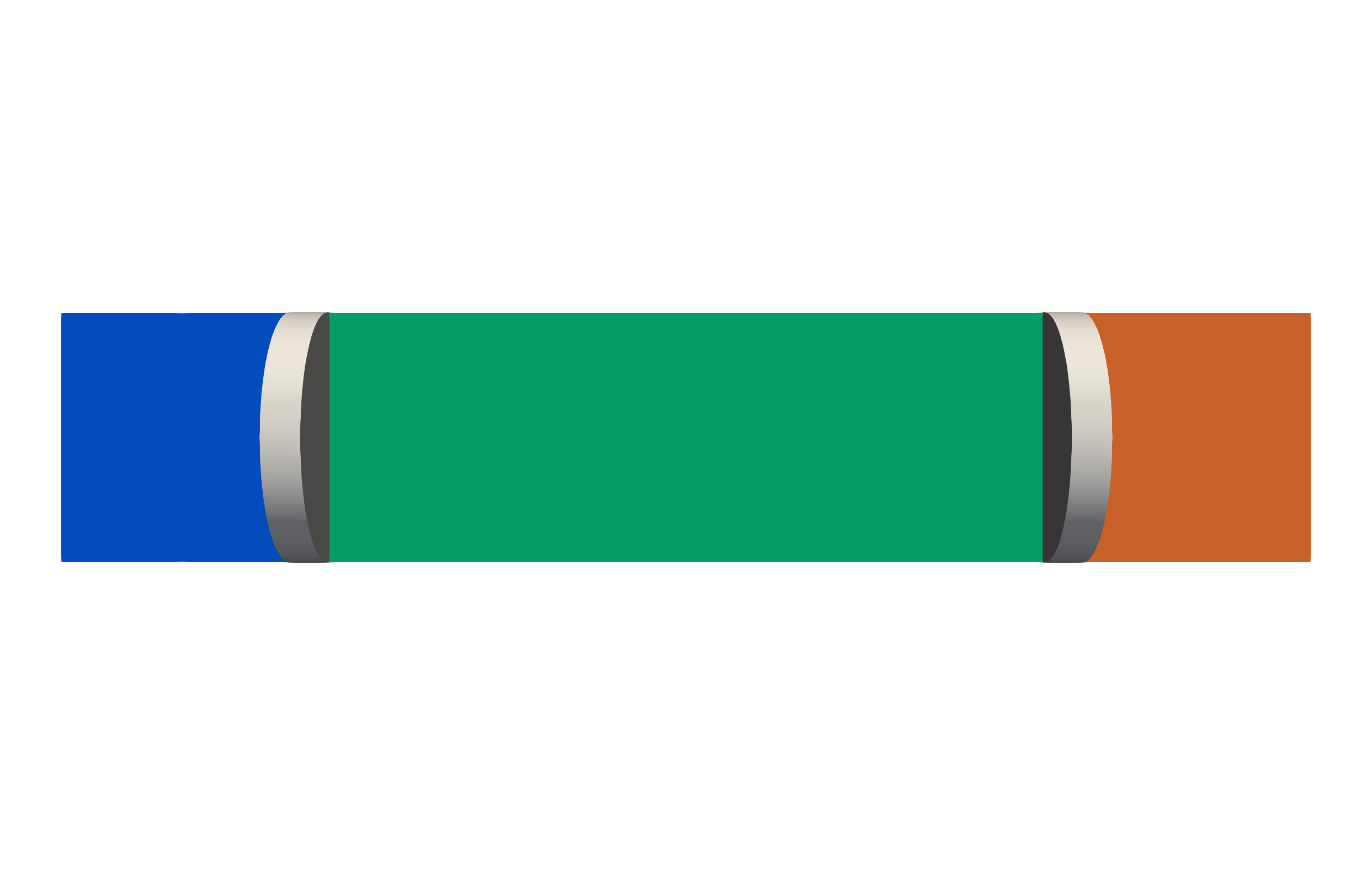}
        \caption{$t = \SI{0.02}{\second}$}
    \end{subfigure}
    
    \vspace{1em}
    \begin{subfigure}{\textwidth}
        \includegraphics[width=0.49\textwidth,
                         trim={0in 7in 0in 7in},
                         clip]{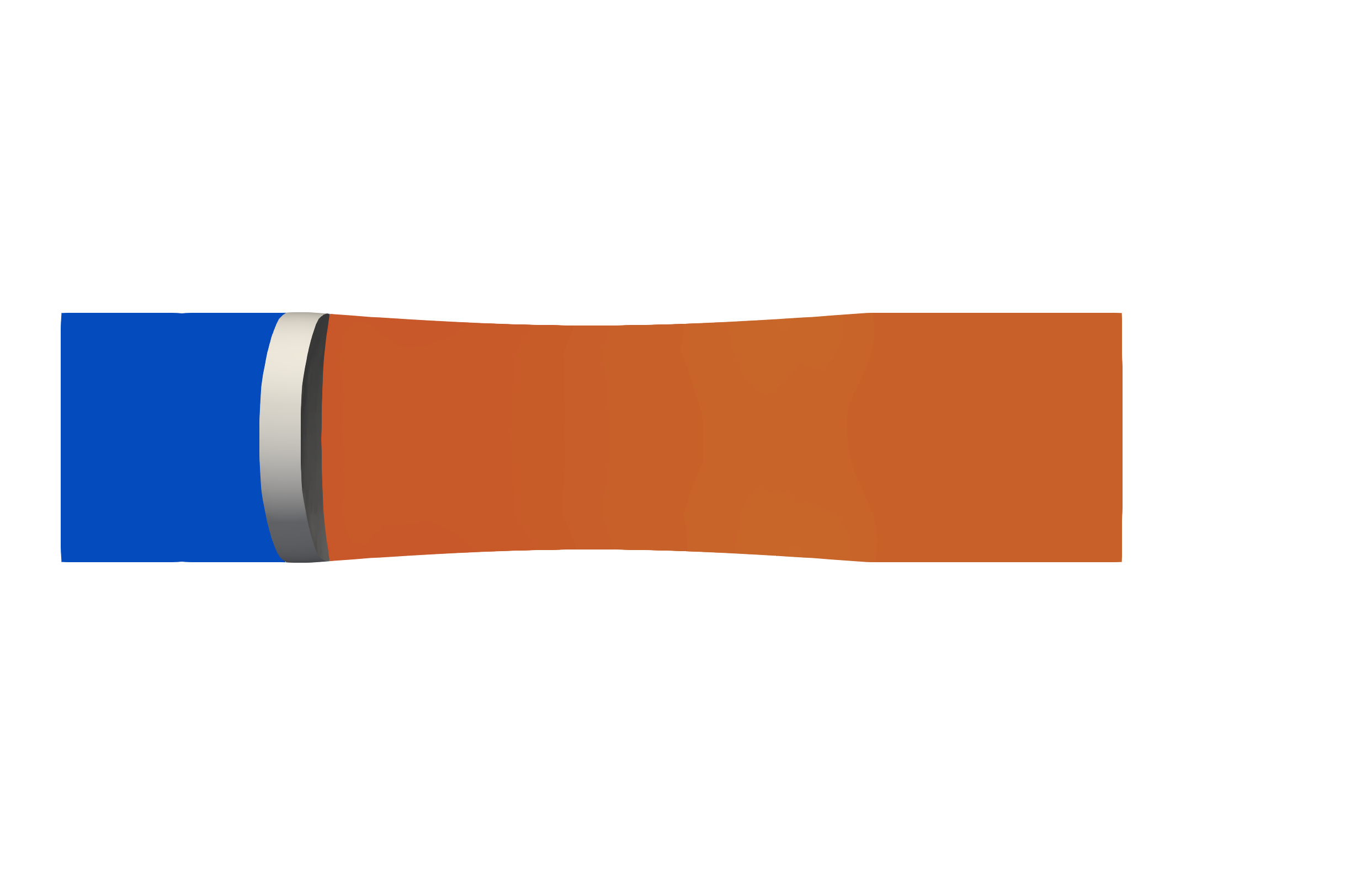}
        \includegraphics[width=0.49\textwidth,
                         trim={0in 7in 0in 7in},
                         clip]{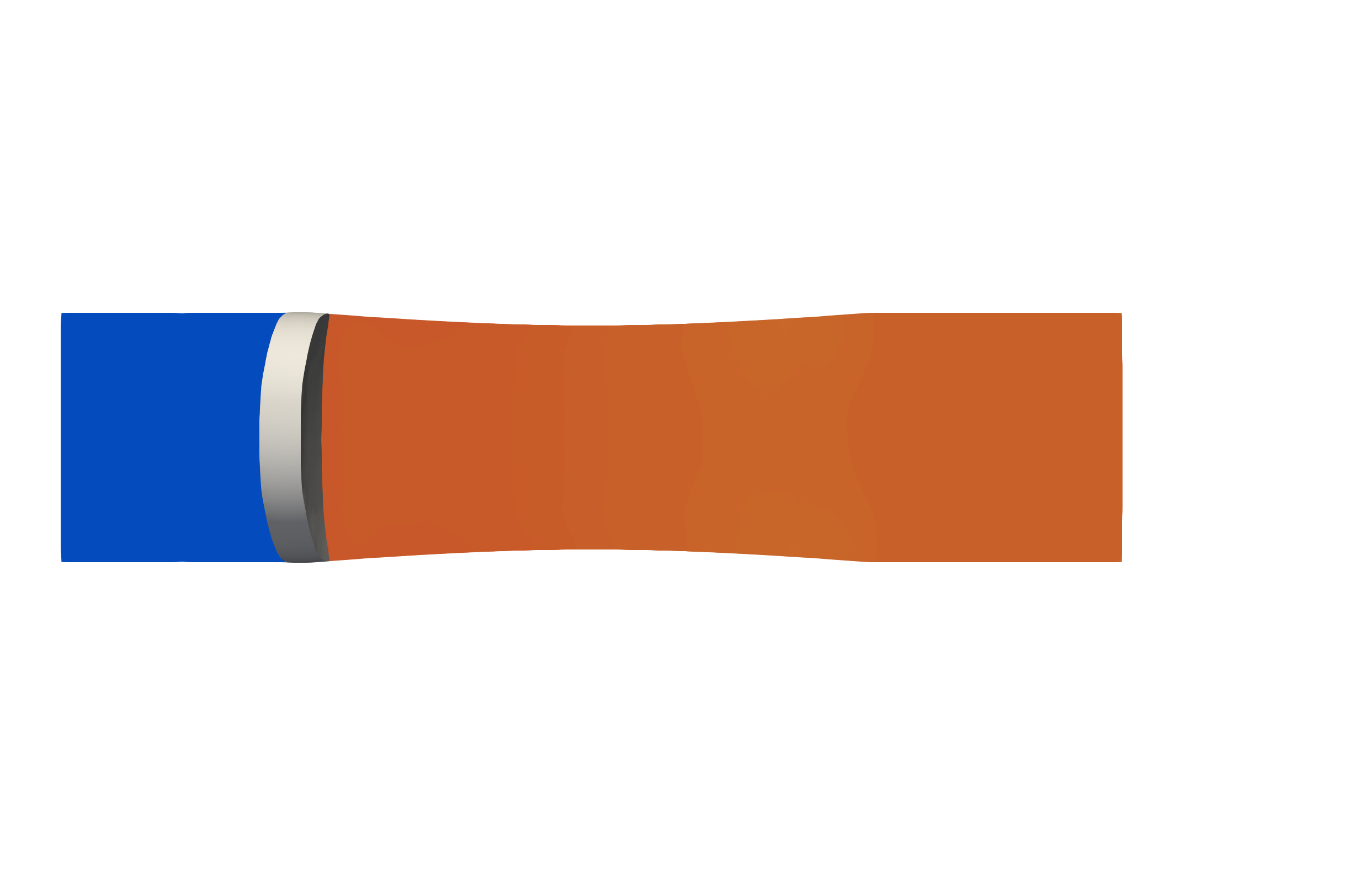}
        \caption{$t = \SI{0.08}{\second}$}
    \end{subfigure}
    
    \vspace{1em}    
    \begin{subfigure}{\textwidth}
        \includegraphics[width=0.49\textwidth,
                         trim={0in 7in 0in 7in},
                         clip]{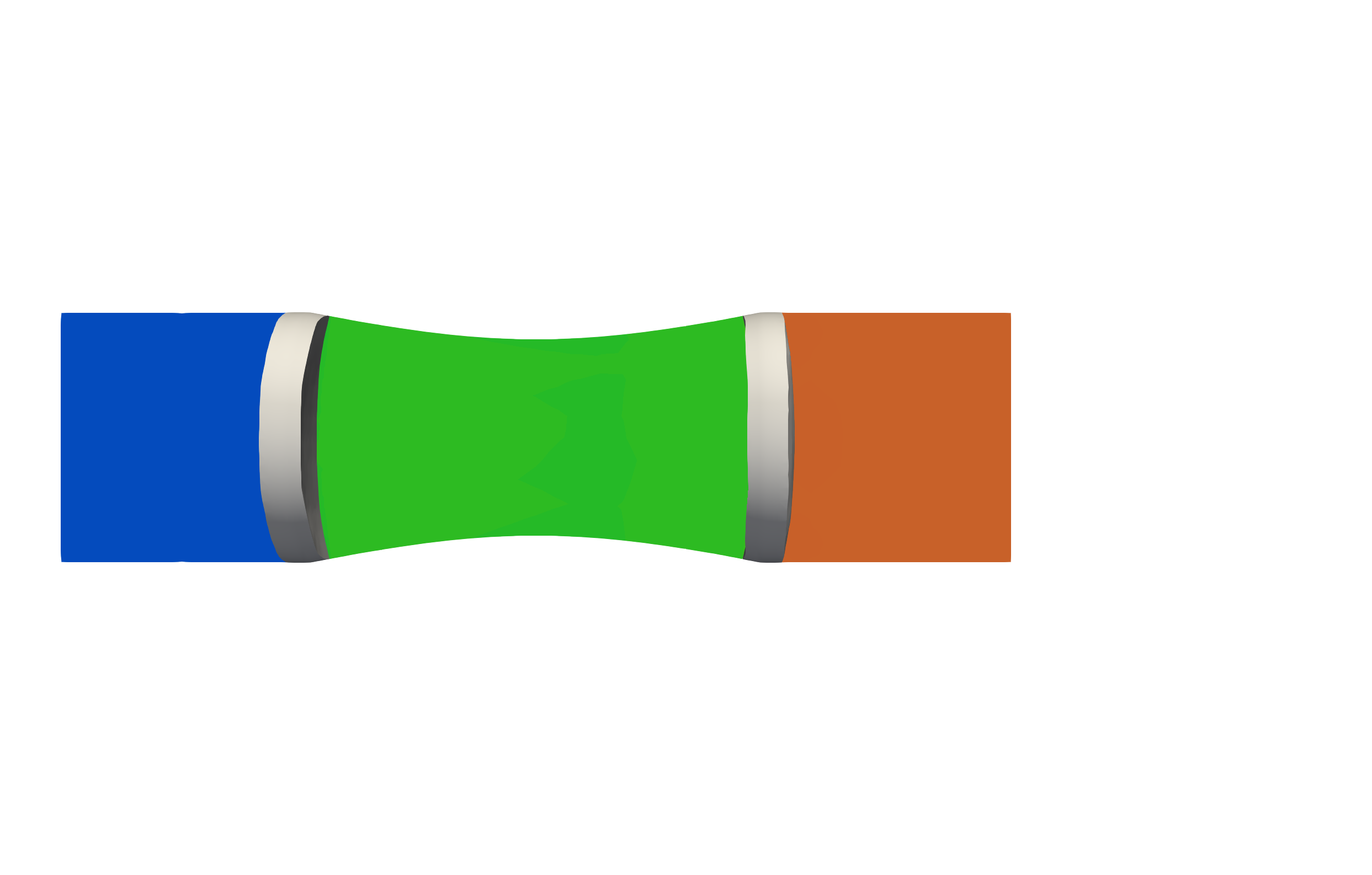}
        \includegraphics[width=0.49\textwidth,
                         trim={0in 7in 0in 7in},
                         clip]{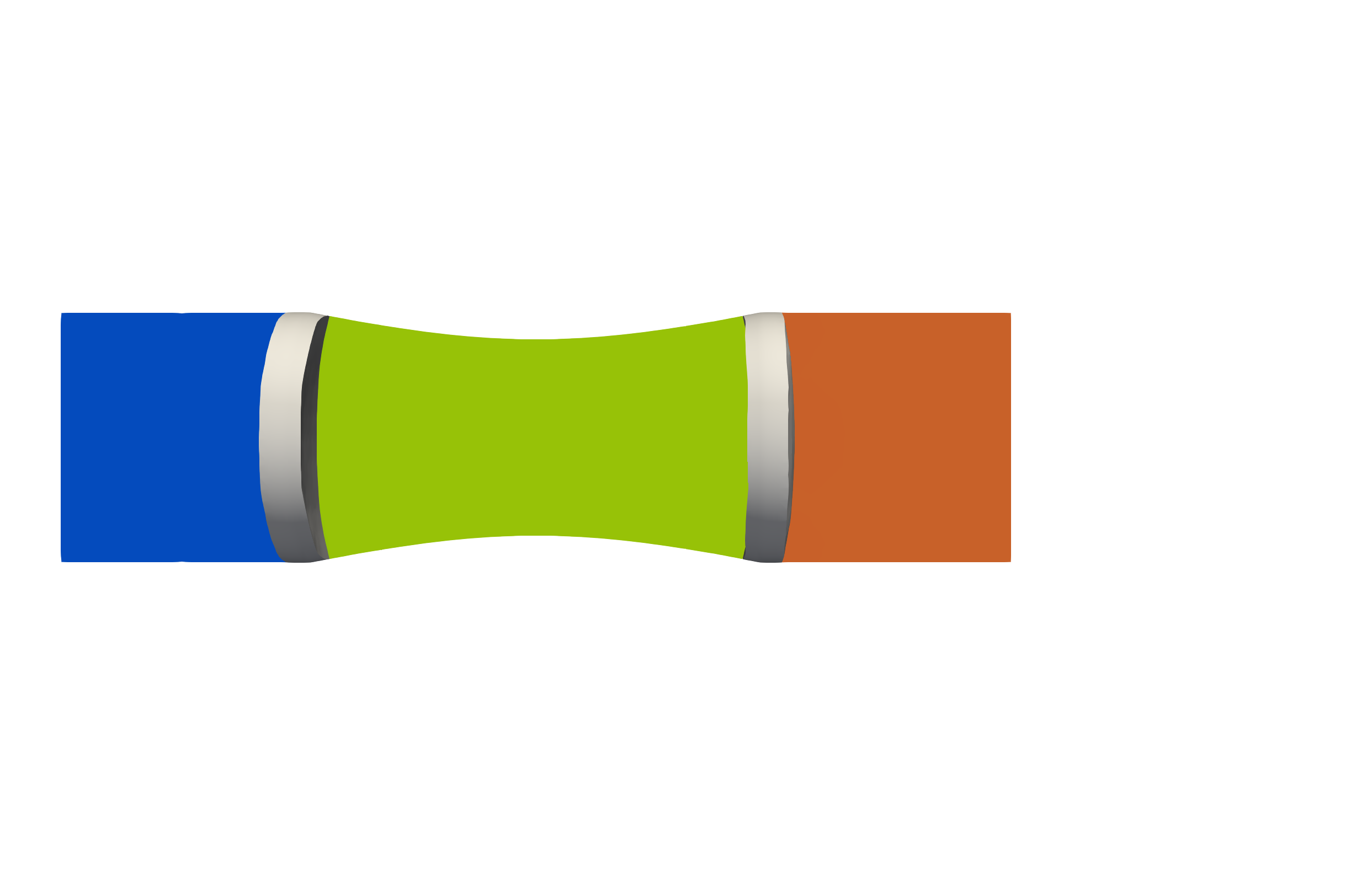}
        \caption{$t = \SI{0.12}{\second}$}
    \end{subfigure}
    
    \vspace{1em}
    \begin{subfigure}{\textwidth}
        \includegraphics[width=0.49\textwidth,
                         trim={0in 7in 0in 7in},
                         clip]{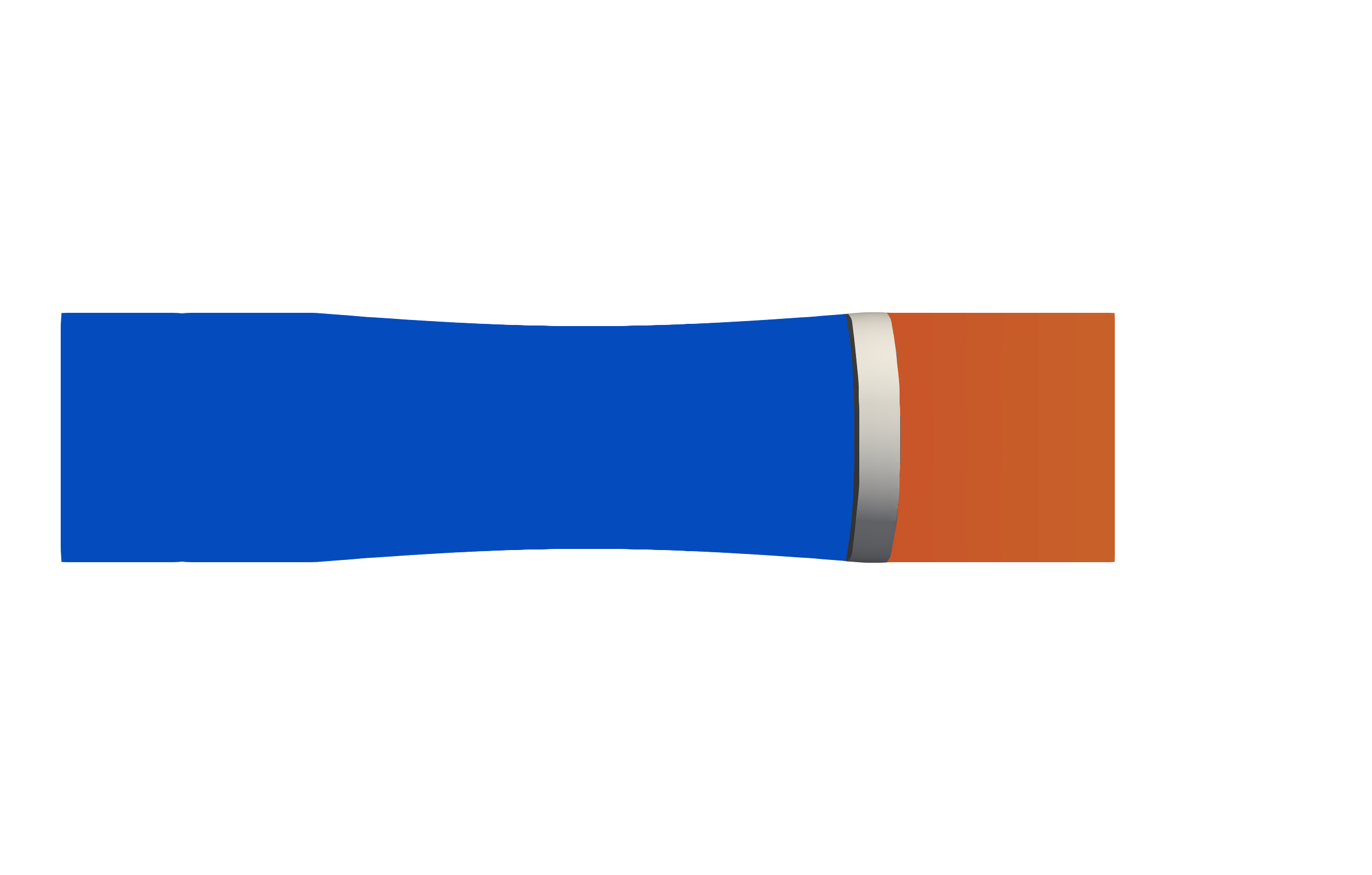}
        \includegraphics[width=0.49\textwidth,
                         trim={0in 7in 0in 7in},
                         clip]{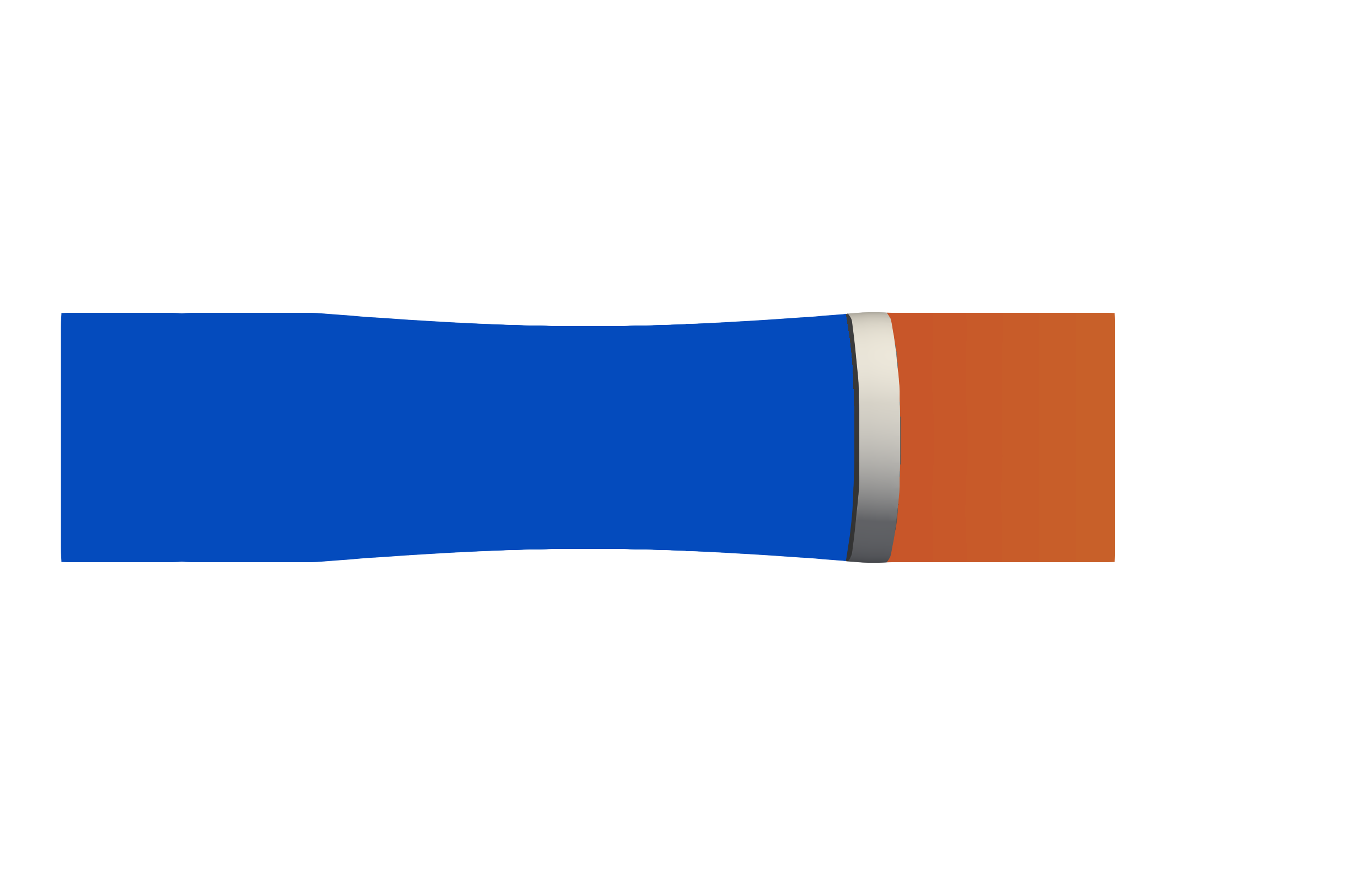}
        \caption{$t = \SI{0.17}{\second}$}
    \end{subfigure}
    
    \vspace{1em}
    \begin{subfigure}{0.25\textwidth}
        \includegraphics[width=\textwidth]{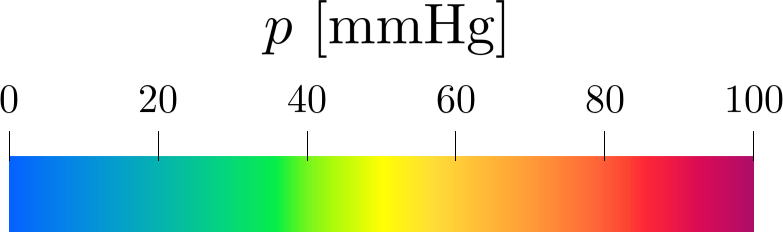}
    \end{subfigure}
    
    \caption{Test B. Snapshots of the pressure over one longitudinal slice of the domain, simulated using the RIIS (left) and ARIIS (right) method. The snapshots are taken at the midpoint of isovolumetric contraction (a), ejection (b), isovolumetric relaxation (c) and filling (d). The domain is warped according to the displacement $\mathbf d$ defined in \eqref{eq:displacement-cardiocylinder}.}
    \label{fig:cardiocylinder-solution-pressure}
\end{figure}

\begin{figure}
    \centering
    \begin{subfigure}{0.49\textwidth}
        \centering
        \textbf{RIIS}
    \end{subfigure}
    \begin{subfigure}{0.49\textwidth}
        \centering
        \textbf{ARIIS}
    \end{subfigure}
    
    \vspace{1em}
    \begin{subfigure}{\textwidth}
        \includegraphics[width=0.49\textwidth]{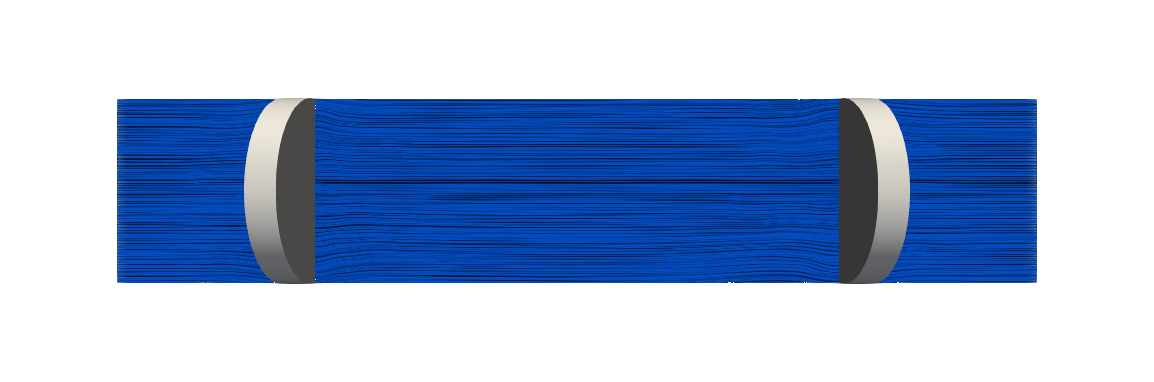}
        \includegraphics[width=0.49\textwidth]{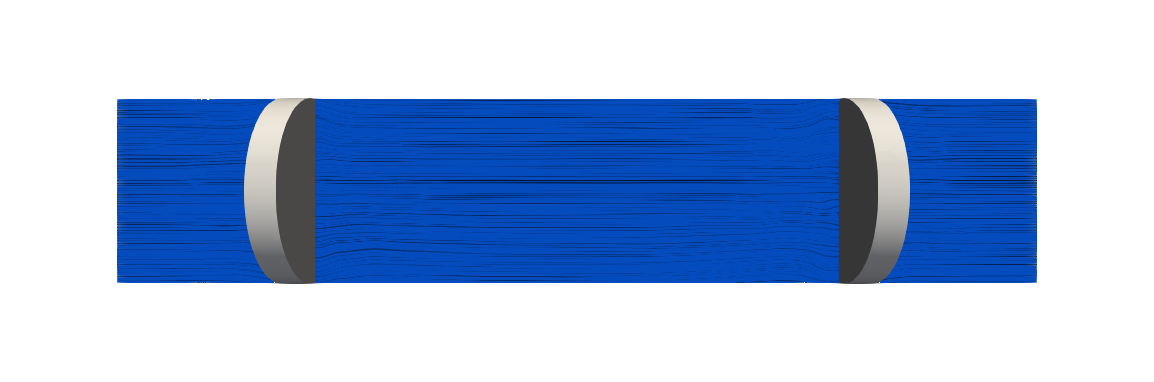}
        \caption{$t = \SI{0.02}{\second}$}
    \end{subfigure}
    
    \vspace{1em}
    \begin{subfigure}{\textwidth}
        \includegraphics[width=0.49\textwidth]{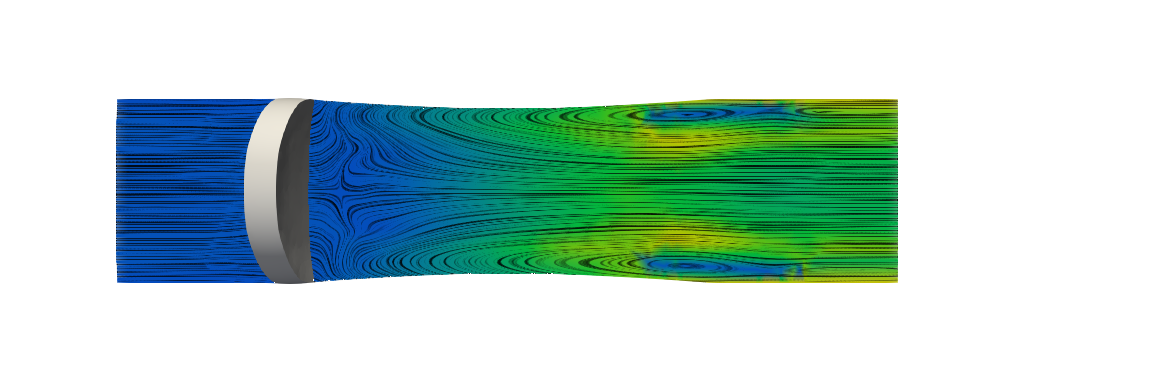}
        \includegraphics[width=0.49\textwidth]{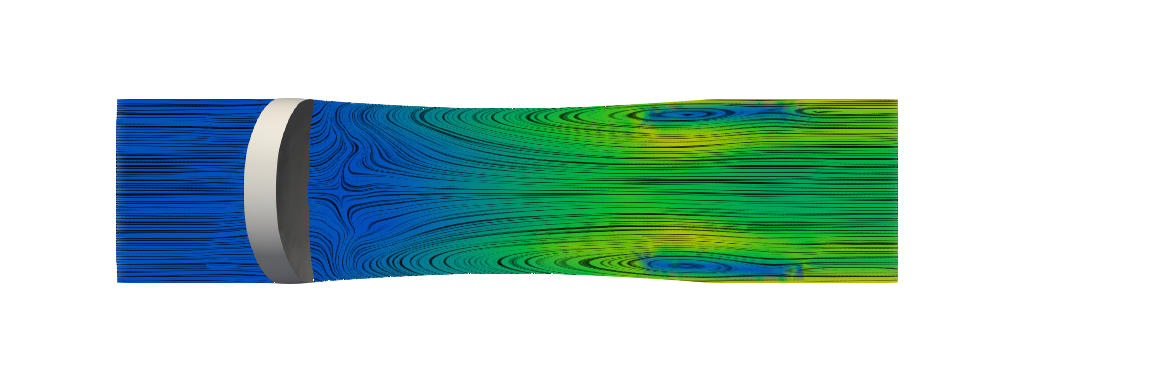}
        \caption{$t = \SI{0.08}{\second}$}
    \end{subfigure}
    
    \vspace{1em}
    \begin{subfigure}{\textwidth}
        \includegraphics[width=0.49\textwidth]{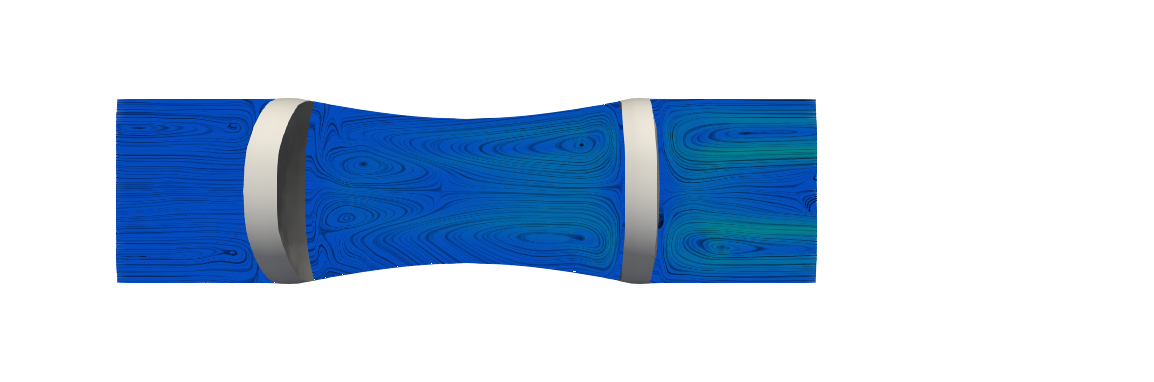}
        \includegraphics[width=0.49\textwidth]{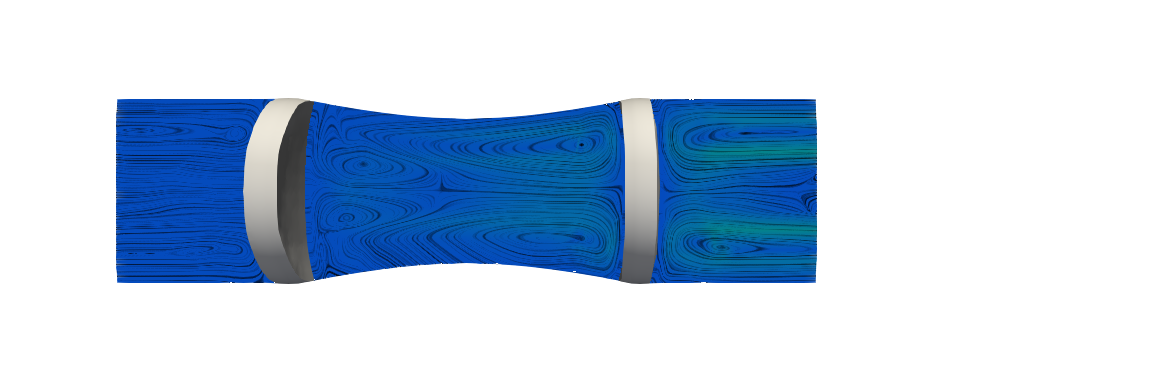}
        \caption{$t = \SI{0.12}{\second}$}
    \end{subfigure}
    
    \vspace{1em}
    \begin{subfigure}{\textwidth}
        \includegraphics[width=0.49\textwidth]{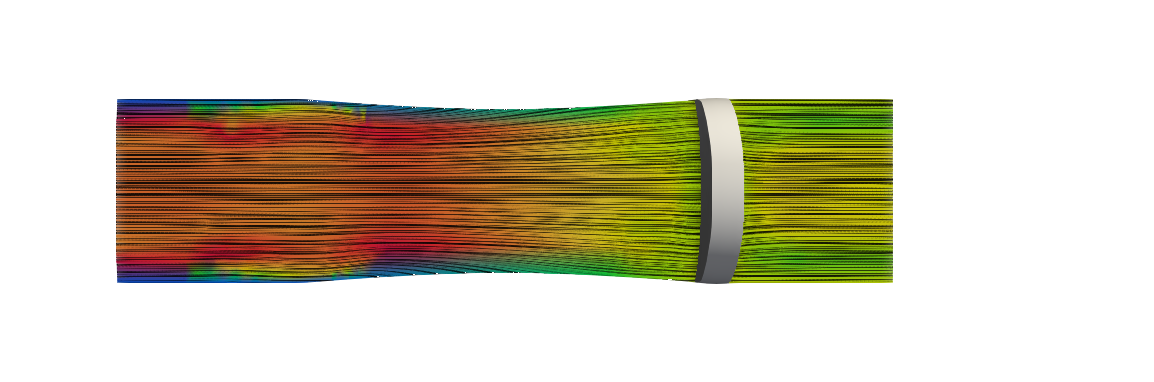}
        \includegraphics[width=0.49\textwidth]{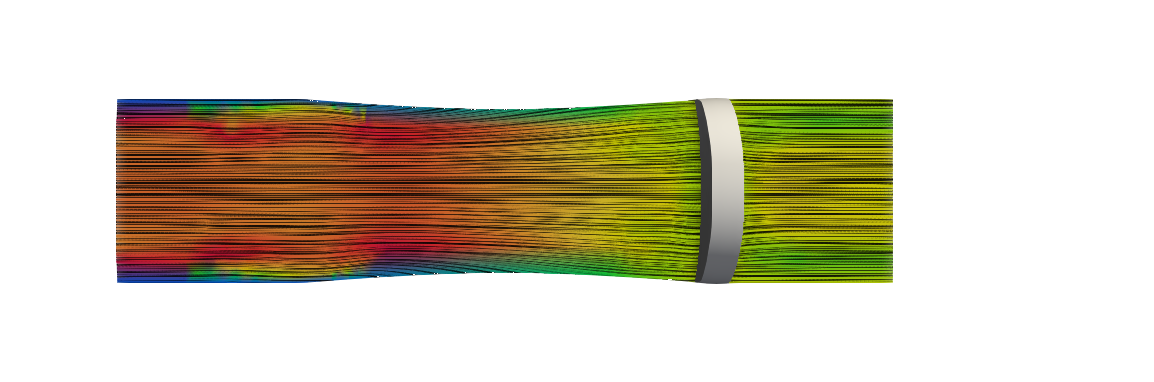}
        \caption{$t = \SI{0.17}{\second}$}
    \end{subfigure}
    
    \vspace{1em}
    \begin{subfigure}{0.25\textwidth}
        \includegraphics[width=\textwidth]{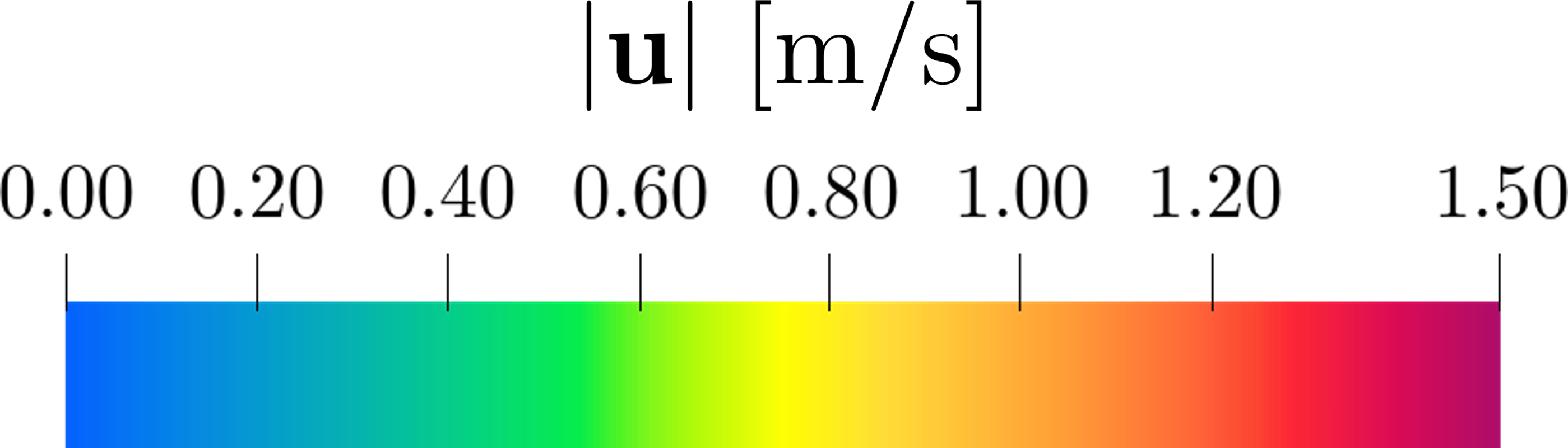}
    \end{subfigure}
    
    \caption{Test B. Snapshots of the velocity magnitude over one longitudinal slice of the domain, simulated using the RIIS (left) and ARIIS (right) methods. The snapshots are taken at the midpoint of isovolumetric contraction (a), ejection (b), isovolumetric relaxation (c) and filling (d). The domain is warped by the displacement $\mathbf d$ defined in \eqref{eq:displacement-cardiocylinder}, and the velocity magnitude is superimposed with the surface LIC rendering of the flow field.}
    \label{fig:cardiocylinder-solution-velocity}
\end{figure}

\begin{figure}
    \centering
    \includegraphics[width=0.7\textwidth]{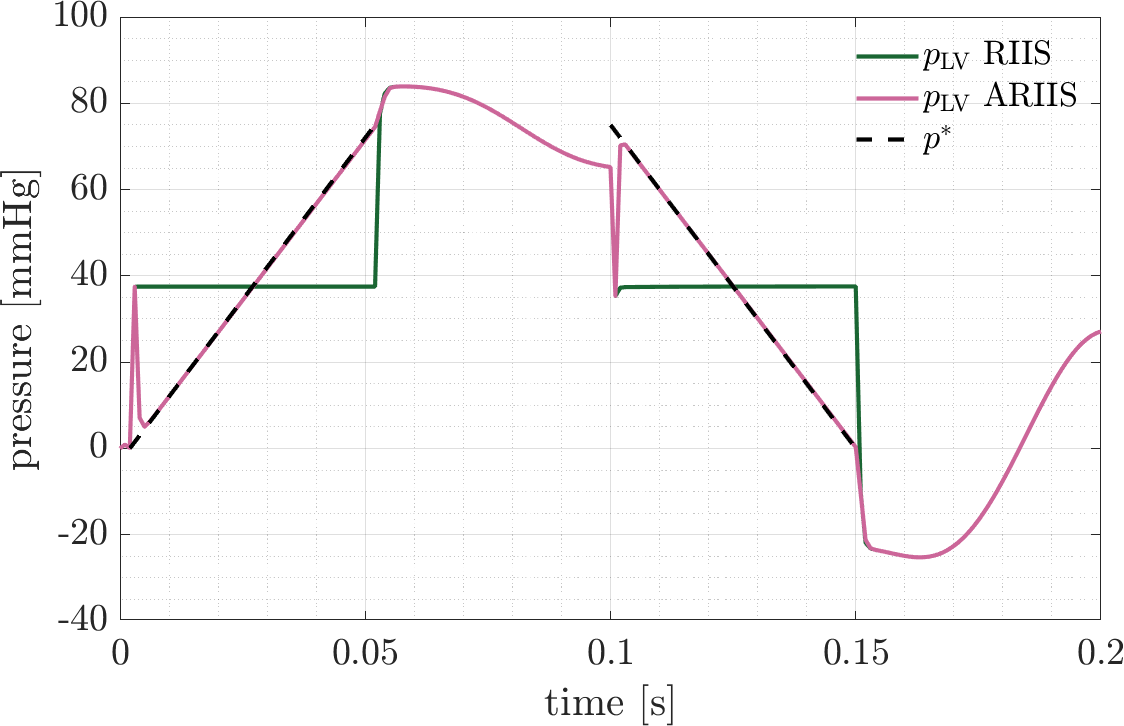}
    \caption{Test B. Ventricular and reference pressures, with RIIS and ARIIS methods, using resistance $R = \SI{e4}{\kilo\gram\per\metre\per\second}$ and $\varepsilon = \SI{0.002}{\metre}$.}
    \label{fig:cardiocylinder-riis-vs-ariis}
\end{figure}

% RIIS vs. ARIIS.
\Cref{fig:cardiocylinder-solution-pressure,fig:cardiocylinder-solution-velocity} report snapshots of pressure and velocity in the solution, computed using RIIS and ARIIS. We can observe that the two methods yield equivalent results outside the isovolumetric phases. Differently, when both valves are closed, a considerably different pressure can be observed.
Similar conclusions can be drawn from the plots reported in \cref{fig:cardiocylinder-riis-vs-ariis}, representing the average ventricular pressure over time for Test B, using RIIS and ARIIS, setting $R = \SI{e4}{\kilo\gram\per\metre\per\second}$ and $\varepsilon = \SI{0.002}{\metre}$. The ARIIS simulation yields a pressure that closely follows the provided reference pressure $p^*$ during isovolumetric phases. Conversely, outside the isovolumetric phases, the two methods correctly produce the same result.

\begin{figure}
    \centering
        \includegraphics[width=0.4\textwidth]{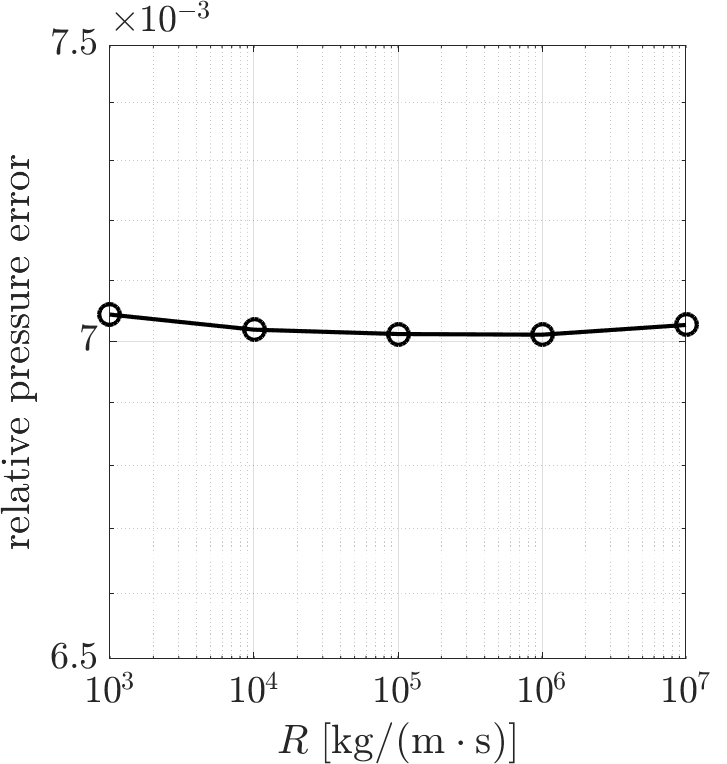}
    \caption{Relative error between the computed pressure during isovolumetric phases and the reference pressure $p^*$, with varying resistance $R$, $\varepsilon = \SI{0.002}{\metre}$ and minimum mesh size $h_\text{min} = \SI{0.001}{\metre}$.}            \label{fig:cardiocylinder-ariis-resistance}
\end{figure}

% Varying the resistance.
We carry out numerical simulations with the ARIIS method by varying the resistance coefficient $R$ over several orders of magnitude and computing the relative pressure error \eqref{eq:relative-pressure-error} during isovolumetric phases. We report the results in \cref{fig:cardiocylinder-ariis-resistance}. As before, we observe that the reference pressure is matched accurately during isovolumetric phases, regardless of the value of $R$.

\FloatBarrier

\subsection{Test C: application to a cardiac case}
\label{sec:cfd-lh}

\begin{figure}[t]
    \centering
    \begin{subfigure}{0.45\textwidth}
        \includegraphics[width=\textwidth]{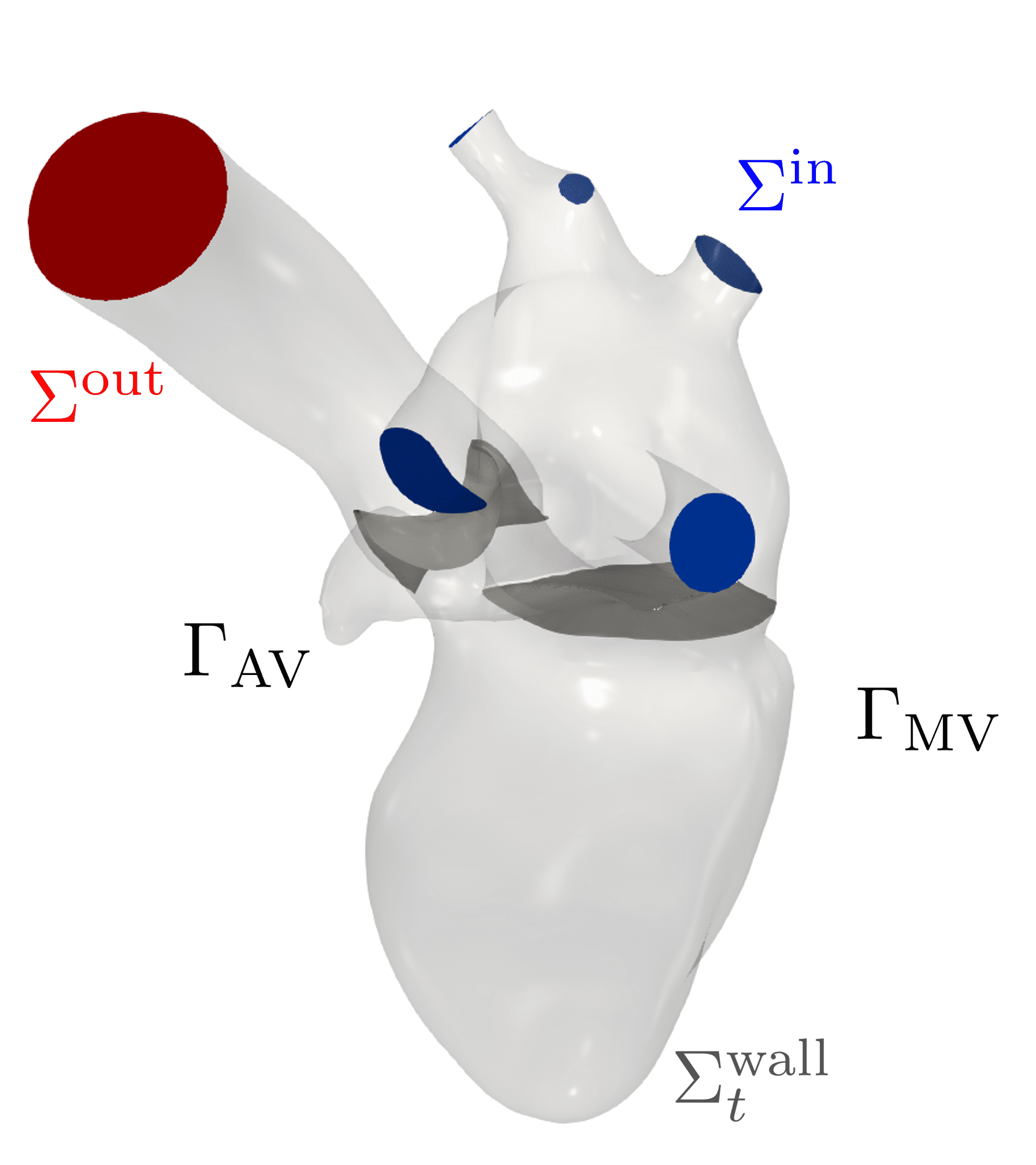}
        \caption{}
        \label{fig:lh-cfd-domain}
    \end{subfigure}
    \begin{subfigure}{0.45\textwidth}
        \includegraphics[width=\textwidth]{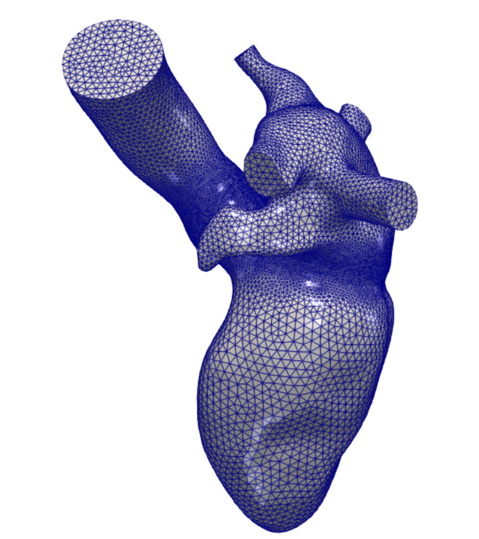}
        \caption{}
        \label{fig:lh-cfd-mesh}
    \end{subfigure}
    \caption{Test C. Left heart CFD domain highlighting boundary portions and immersed surfaces in their closed configurations (a); tetrahedral mesh generated for the CFD simulation (b).}
    \label{fig:lh-cfd-domain-mesh}
\end{figure}

In this section, we apply the ARIIS method to a realistic cardiac case. We use the CFD model of a healthy left heart developed in \cite{zingaro2022geometric}. It consists of the 3D fluid dynamics model \eqref{eq:ns-ariis} coupled to the surrounding circulation (described by a 0D closed-loop model \cite{blanco2010computational, hirschvogel2017monolithic, REGAZZONI2022111083}) and driven by a cardiac electromechanical model \cite{REGAZZONI2022111083}.

\begin{table}[t]
\centering
	\begin{tabular}{c|S S S S S}
% 	\hline
	$\mathrm k$ & 
	$R_{\mathrm k}$ [\si{\kilo\gram\per\metre\per\second}] & 
	$\varepsilon_\mathrm k$ [\si{\milli\metre}] &
	$|\Gamma_{\mathrm k}|$ [\si{\cm\squared}] &
	\text{clos. time} [\si{\second}] &
	\text{open. time} [\si{\second}]
	\\
	\hline
	MV & 1e4 & 1.0 & 12.11 & 0.04725 & 0.49350 \\
	AV & 1e4 & 1.0 & 5.41 & 0.38850 & 0.10600\\
% 	\hline
	\end{tabular}
	\caption{Test C. Parameters of the RIIS and ARIIS methods in the left heart CFD simulations.}
	\label{tab:lh-riis-ariis}
\end{table}
\begin{table}[t!]
	\centering
\begin{tabular}{ccc|c|ccc}
% \hline
% 	\multicolumn{3}{|c|}{$h$} & {cells} & \multicolumn{3}{c|}{{DOFs ($ \mathbb{P}_1 - \mathbb{P}_1 $)}}\\
% 	\multicolumn{3}{|c|}{\footnotesize{[\si{\milli\metre}]}} & [-] &  \multicolumn{3}{c|}{[-]} \\
	\multicolumn{3}{c|}{$h$ \footnotesize{[\si{\milli\metre}]}} & {cells \footnotesize{[-]}} & \multicolumn{3}{c}{{DOFs ($ \mathbb{P}_1 - \mathbb{P}_1 $)} \footnotesize{[-]}}\\
	\footnotesize{min} & \footnotesize{avg} & \footnotesize{max}  && \footnotesize{$\mathbf u$} & \footnotesize{$p$} & \footnotesize{total}\\
% 	\hline
	\hline
	\num{0.44} & \num{1.44} & \num{6.40} & \num{645699} & \num{333243} & \num{111081} & \num{444324} \\
% 	\hline
\end{tabular}
	\caption{Test C. Mesh details for the left heart CFD simulations.}
	\label{tab:lh-mesh}
\end{table}

\begin{figure}[t]
    \centering
    \includegraphics[width=\textwidth]{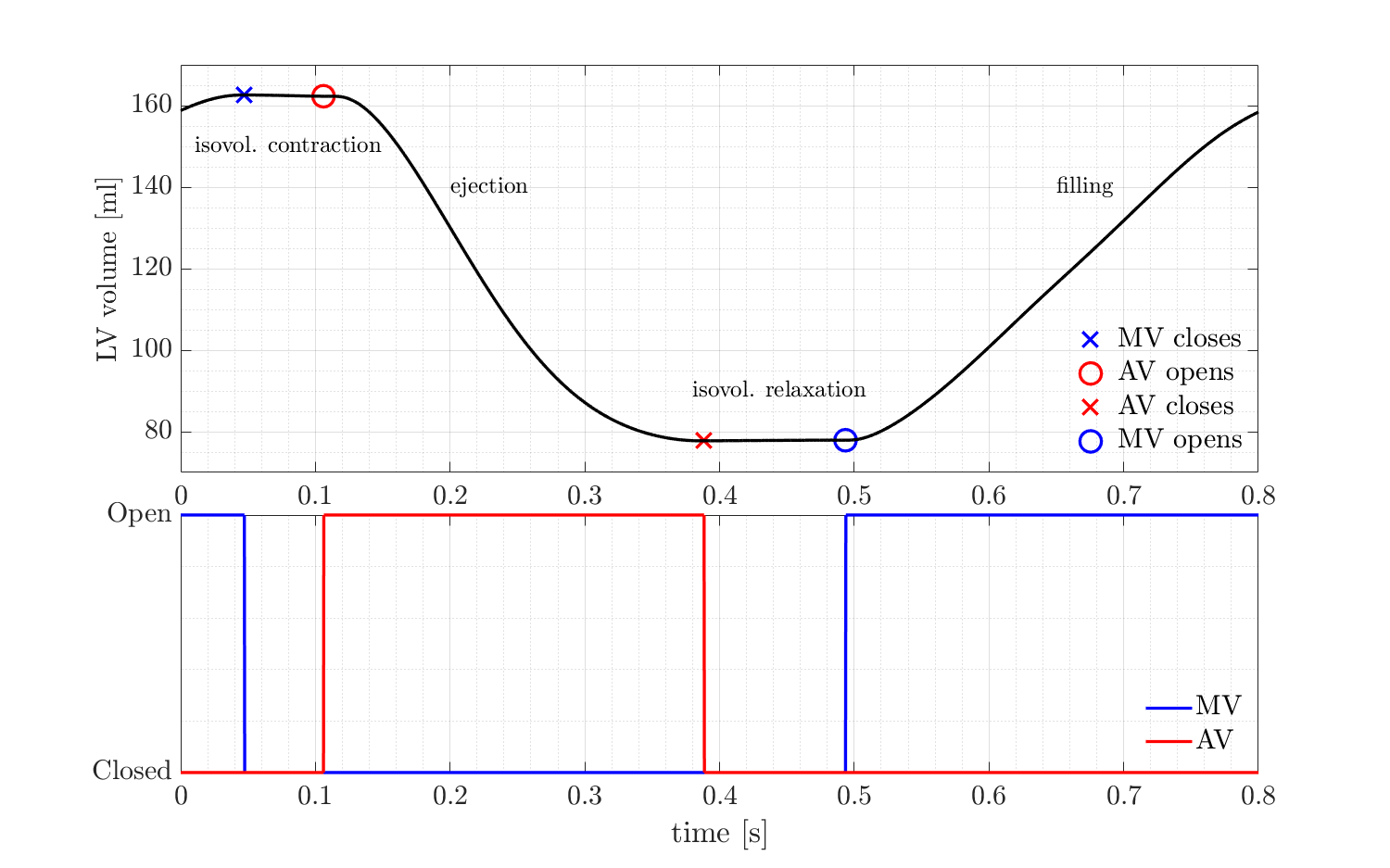}
    \caption{Test C. Volume of left ventricle, with opening and closing times for valves, and valve states.}
    \label{fig:lh-lv-volume}
\end{figure}
We consider a realistic left heart geometry provided by Zygote \cite{zygote}, representing an accurate 3D model of the heart obtained with CT scan data. We report the domain in \Cref{fig:lh-cfd-domain}: its boundary is split as $\partial \Omega_t = \Sigma^\mathrm{in} \cup \Sigma^\mathrm{out} \cup \Sigma^\mathrm{wall}_t $, where $\Sigma^\mathrm{in}$ is the set of pulmonary veins inlet sections, $\Sigma^\mathrm{out}$ the outlet section of the ascending aorta and $\Sigma^\mathrm{w}_t $ the wall (endocardium). In addition, we display the immersed surfaces $\Gamma_\mathrm{MV}$ and $\Gamma_\mathrm{AV}$ in their closed configurations.  

We set Neumann boundary conditions on the inlet and outlet sections of the domain by prescribing the pressure coming from the coupling between the 3D and the 0D circulation model, as explained in \cite{zingaro2021hemodynamics}. To prescribe the displacement field on the endocardium of the LV, we carry out an electromechanical simulation with the ventricular model proposed in \cite{REGAZZONI2022111083}. We report the complete setup of the electromechanical model in \Cref{appendix:em}. Moreover, since the focus of the paper is the correct estimation of the ventricular pressure only, we neglect the motion of the remaining part of the domain by setting homogeneous Dirichlet boundary conditions on the wall of the left atrium and the ascending aorta. To avoid mesh element distortion, for the computation of $K$ in \cref{eq:lifting}, we use the boundary-based stiffening approach proposed in \cite{jasak2006automatic}. 

We generate the tethraedral mesh of the left heart displayed in \Cref{fig:lh-cfd-mesh} with \texttt{vmtk} \cite{vmtk} using the methods and tools discussed in \cite{fedele2021polygonal,zingaro2022geometric}. Mesh details are summarized in \Cref{tab:lh-mesh}. We use as time-step size $\Delta t = \SI{2.5e-4}{\second}$. Since the electromechanical simulation has a much larger timestep than the CFD one, we use smoothing splines \cite{de1978practical} to approximate the electromechanical displacement field in time. 

The values of $R_\mathrm{k}$ and $\varepsilon_\mathrm{k}$ of the RIIS method are provided in \Cref{tab:lh-riis-ariis}. 
These values of $\varepsilon_\mathrm{k}$ and $R_\text{k}$ prevent flow through the closed immersed surfaces \cite{fedele2017patient}. Moreover, following \cite{fedele2017patient}, we choose $\varepsilon_\k$ to guarantee that $\varepsilon_\k \geq 1.5\, h_\text{min}$, where $h_\text{min}$ is the minimum mesh size in the valves region. 
Since the condition number of the linear system associated to the FE discretization of \eqref{eq:ns-ariis} increases as the ratio $R_\mathrm{k}/\varepsilon_\mathrm{k}$ increases, we choose the minimum value of $R_{\mathrm k}$ that guarantees impervious valves, as in \cite{zingaro2022geometric}. In \Cref{tab:lh-riis-ariis}, we also report the areas of the valve sections needed for the ARIIS method. Moreover, as reference pressure $p^*(t)$, we use the one computed in the 3D-0D electromechanical ventricular model \cite{REGAZZONI2022111083}.

Numerical simulations are run in parallel using 48 cores from the GALILEO100 supercomputer at the CINECA supercomputing center.

\subsubsection{Comparison of RIIS and ARIIS methods}

We carry out numerical simulations with the RIIS and the ARIIS methods. We simulate a single heartbeat of period $T=0.8$ s. In \Cref{fig:lh-lv-volume}, we display the LV volume with the four heartbeat phases, along with the times corresponding to the begin and end of isovolumetric phases. We open and close the valves instantaneously (i.e. in one time step) at these times, also reported  in \Cref{tab:lh-riis-ariis}.  {As reference pressure ($p^*$) for the ARIIS method, we use the LV pressure coming from the 3D cardiac electromechanical simulation coupled to the 0D cardiocirculatory model \cite{REGAZZONI2022111083}}.

\begin{figure}
  \centering
        \includegraphics[width=0.9\textwidth]{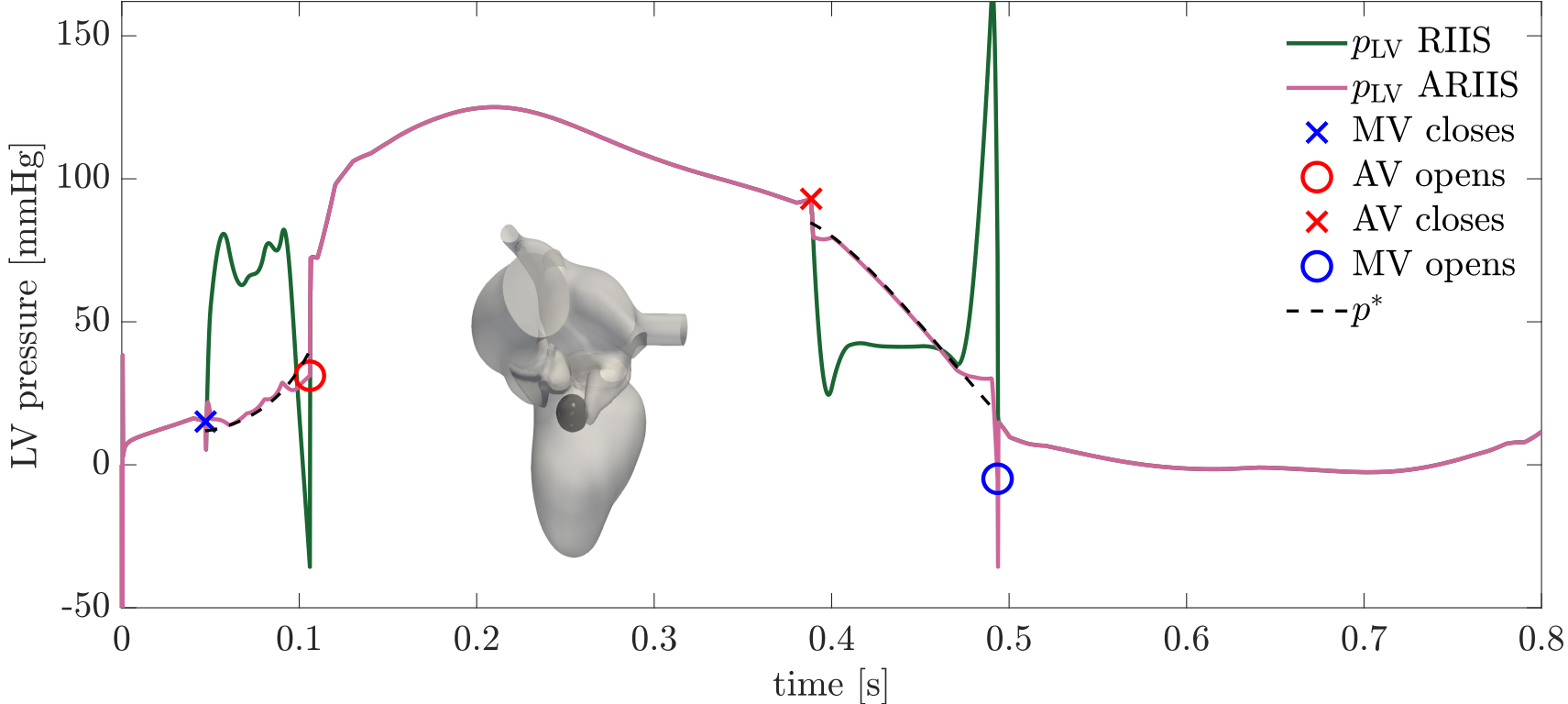}
    \caption{Test C. Ventricular and reference pressures for the left heart test case, with RIIS and ARIIS methods. $p_\mathrm{LV}$ is computed space-averaging the fluid pressure in the black control volume in the left ventricle.}
    \label{fig:lh-cfd-pressure-riis-ariis}
\end{figure}

We display the ventricular pressure with the RIIS and ARIIS methods in \Cref{fig:lh-cfd-pressure-riis-ariis}. We compute it by space-averaging the pressure in a control volume downwind of the MV. The RIIS method is not able to correctly capture the left ventricular pressure, yielding arbitrary pressure values during the isovolumetric phases, with unphysical oscillations. Differently, with the correction term introduced by the ARIIS method, the ventricular pressure follows the expected trend given by $p^*$. In addition, out of the isovolumetric phases, the pressure fields are almost identical between RIIS and ARIIS methods. Indeed, the correction term is active in the isovolumetric phases only, and it does not influence the remaining phases of the heart cycle, yielding a maximum discrepancy of $0.23$ mmHg.

\begin{figure}
    \centering
    \includegraphics[ width=\textwidth]{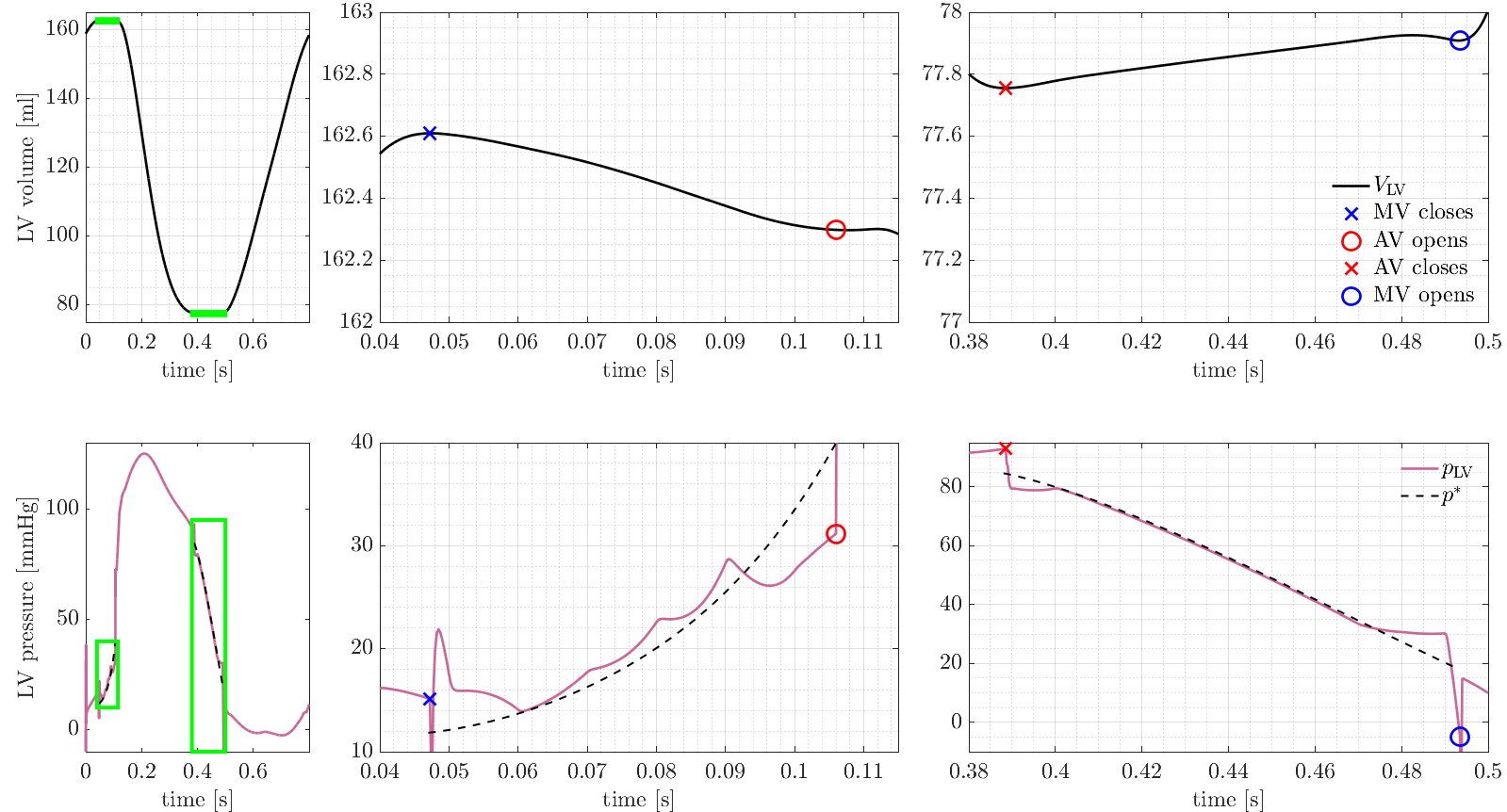}
    \caption{Test C. Ventricular volume and pressures (obtained through the ARIIS method), with zoom on the isovolumetric phases.}
    \label{fig:lh-volume-pressure-zooms}
\end{figure}

Furthermore, as shown in \Cref{fig:lh-volume-pressure-zooms}, the largest discrepancies between $p_\text{LV}$ and $p^*$ in the ARIIS case is attained at the end of the isovolumetric phases. These discrepancies are related to the fact that the isovolumetric phases in realistic cardiac simulations are not exactly volume preserving. This happens due to, on the one hand, the projection of the displacement from the electromechanics (or imaging data) onto the fluid dynamics mesh and, on the other hand, the lifting problem in \eqref{eq:lifting} that does not guarantee, a priori, any kind of volume conservation in the LV subdomain. Moreover, the displacement is characterized by small oscillations in time -- introduced by the smoothing splines -- that yield oscillations in the ventricular volume as well. 
Nonetheless, differently from the standard RIIS method, the proposed augmented approach allows to simulate the isovolumetric phases, with a pressure evolution that is much more similar to the heart physiology. 

\begin{figure}
  \centering
        \includegraphics[width=0.8\textwidth]{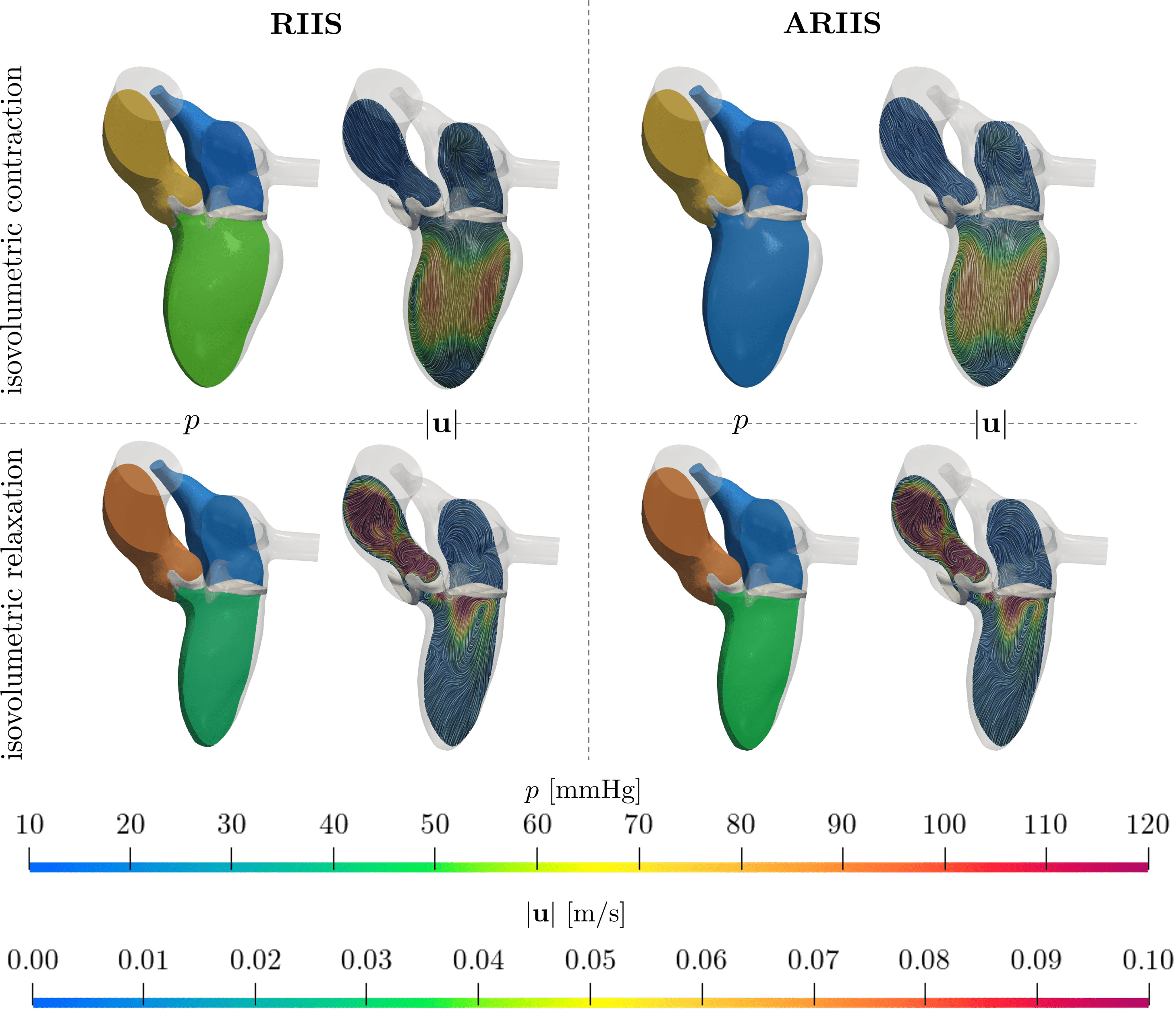}
    \caption{Test C. Comparison between RIIS and ARIIS methods during isovolumetric phases: pressure on a clip in the LV apico-basal direction and a section colored according to velocity magnitude with a surface LIC representation.}
    \label{fig:lh-clip-pressure-velocityLIC-riis-ariis}
\end{figure}

In \Cref{fig:lh-clip-pressure-velocityLIC-riis-ariis}, we show the pressure field (in mmHg) on a clip in the LV apico-basal direction during the isovolumetric phases. The RIIS and ARIIS methods are characterized by different pressures, confirming our previous results. Moreover, we investigate the difference among the two solutions also in terms of velocity field, by showing a surface LIC representation on a slice in the LV apico-basal direction colored with velocity magnitude. Consistently with the findings of \cite{this2020augmented}, we notice that the augmented approach does not impact the velocity field and both solutions reproduce the same flow patterns. More quantitatively, we compute the velocity magnitude in a control volume in the LV. When the augmented formulation is active, we compute a maximum discrepancy between the RIIS and the ARIIS velocities equal to \num{2.21e-04} \si{\metre\per\second}, corresponding to a relative error (divided by the maximum RIIS velocity magnitude) equal to \num{0.29} \%.

\FloatBarrier

\section{Conclusions and limitations}

In this paper, we proposed an augmented version of the Resistive Immersed Implicit Surface (RIIS) method \cite{fedele2017patient} to correctly simulate the heart hemodynamics during isovolumetric phases. This Augmented RIIS (ARIIS) method extends the previously proposed Augmented Resistive Immersed Surface (ARIS) method \cite{this2020augmented} to the case of meshes that are non conforming to cardiac valves. 

Starting from the RIIS method, we derived the correction term required to simulate the intracardiac hemodynamics when both valves are closed. Specifically, we introduced an additional term to the momentum balance of the Navier-Stokes equations that only acts on the valves and is only active during the isovolumetric phases. From the ARIIS derivation, we found that the corrective term depends on the external pressure, the valve areas, the resistive term itself, and a prescribed (reference) pressure representing the intraventricular pressure transient when both valves are closed. The reference pressure can be imposed, for instance, from electromechanical simulations or from patient-specific data.  

We applied the ARIIS method to three different problems: the same cylindrical toy problem introduced in \cite{this2020augmented} for the sake of validation of the proposed method, a novel benchmark problem retaining characteristics of a heart cycle, and the flow in a realistic human left heart geometry (with endocardium displacement obrtained from electromechanical simulations).

All tests showed that the ARIIS method  {yields} a ventricular pressure that closely follows the prescribed reference evolution. Moreover, we found that the accuracy of the results is not  {affected} by resistance coefficient values. On the other hand, we found that the error between measured and prescribed pressure decreases as the ratio $\varepsilon / h_\mathrm{min}$ increases, where $\varepsilon$ is the half-thickness of the valve and $h_\mathrm{min}$ the minimum mesh size. This result is consistent with the RIIS method better capturing the immersed surface as the number of elements in the resistive surface thickness increases.

The ARIIS method is very sensitive to small volume variations and oscillations during isovolumetric phases. Thus, further investigations are advisable for the employment of a better interpolant or approximant (in time) of the input displacement field. Moreover, we observed some mismatch between the fluid pressure and the electromechanical one, yielding an unphysical jump from the isovolumetric contraction to the ejection phase. This mismatch  {suggests} a deeper investigation of the similarities and differences between electromechanics and CFD models, which will be the subject of future work.

To conclude, the standard RIIS method yielded a ventricular pressure  {with large oscillations in time and inconsistent with physiology.} On the contrary, the perturbation term introduced by the proposed ARIIS method provided a valid approach to produce a far more physiological ventricular pressure, and hence to correctly simulate the isovolumetric phases.  

\section*{Acknowledgments}
AZ, LD and AQ received funding from the Italian Ministry of University and Research (MIUR) within the PRIN (Research projects of relevant national interest 2017 “Modeling the heart across the scales: from cardiac cells to the whole organ” Grant Registration number 2017AXL54F). 

MB, IF, LD and AQ acknowledge the ERC Advanced Grant iHEART, ``An Integrated Heart Model for the simulation of the cardiac function'', 2017–2022,  P.I. A. Quarteroni (ERC–2016– ADG, project ID: 740132). 
\begin{center}
    \includegraphics[width=0.7\textwidth]{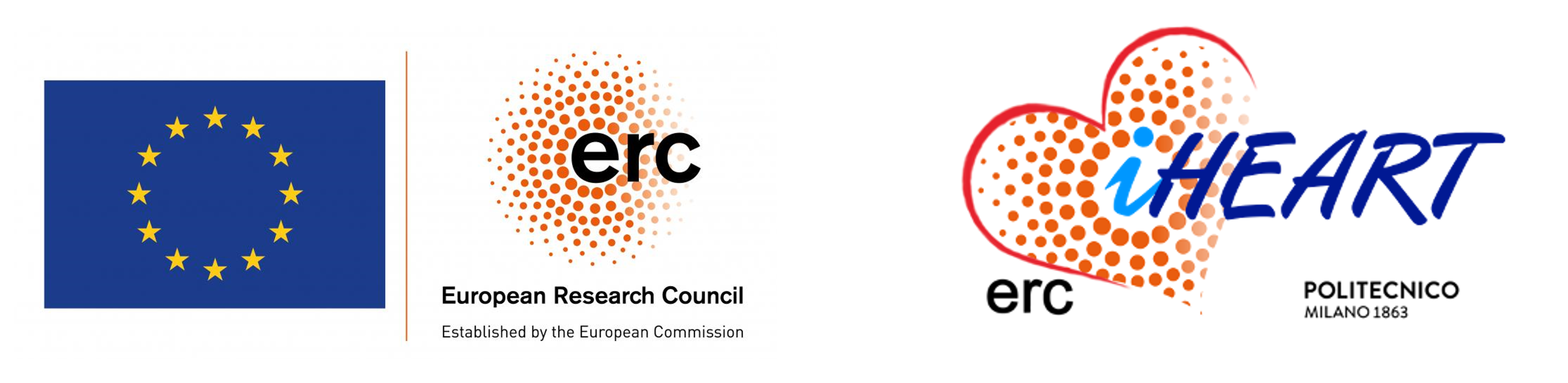}
\end{center}

Finally, we gratefully acknowledge the CINECA award under the ISCRA initiative, for the availability of high performance computing resources and support under the projects IsC87\_MCH, P.I. A. Zingaro, 2021-2022 and IsB25\_MathBeat, P.I. A. Quarteroni, 2021-2022. 

\section*{Ethical statement}
None

\FloatBarrier

\appendix

\section{The electromechanical model}
\label{appendix:em}

We use the electromechanical model developed in \cite{REGAZZONI2022111083}. In the following, we list the parameters employed to carry out the ventricular electromechanical simulation. 

\begin{table}
    \centering

    \begin{tabular}{l| l r S S}
        % \hline
        \textbf{Physics} & \textbf{Parameter} & \multicolumn{2}{c}{\textbf{Value}} 
        \\
        % \hline
        \hline
        \multirow{6}{*}{EP} & \multirow{3}{*}{Conductivities} &
        $\sigma_\mathrm{m}^\mathrm{l}$ & 2.00e-4 & \si{\square\metre\per\second} \\
        & & $\sigma_\mathrm{m}^\mathrm{t}$ & 1.05e-4 & \si{\square\metre\per\second} \\
        & & $\sigma_\mathrm{m}^\mathrm{n}$ & 0.55e-4 & \si{\square\metre\per\second} \\
        \cline{2-5}
        & \multirow{3}{*}{Stimulus} &
        $A_\mathrm{app}$ & 25.71 & \si{\volt\per\second} \\
        & & $\sigma_\mathrm{app}$ & 5e-3 & \si{\metre} \\
        & & $T_\mathrm{app}$ & 3e-3 & \si{\second} \\
        % \hline
        \hline
        \multirow{8}{*}{AFG} &  & $\gamma$ &  30 & \\
        & & $k_d$ &  0.36 & \\
         & & $\alpha_{k_d}$  & -0.2083 & \\
        & &$K_\mathrm{off}$ & 8 & \si{\per\second} \\
        & & $K_\mathrm{basic}$  &  4 &  \si{\per\second} \\
        & & $\mu_{fp}^0$ & 32.255  & \si{\per\second} \\
         & & $\mu_{fp}^1$ & 0.768  & \si{\per\second} \\
         & & $a_\mathrm{XB}$ & 20e8  & \si{\pascal} \\
        \hline
        % \hline
        \multirow{13}{*}{M} & 
        \multirow{8}{3.75cm}{Guccione}
        & $c$ & 8.8e2 & \si{\pascal} \\
        & & $a_\mathrm{ff}$ & 8 & \\
        & & $a_\mathrm{ss}$ & 6 & \\
        & & $a_\mathrm{nn}$ & 3 & \\
        & & $a_\mathrm{fs}$ & 12 & \\
        & & $a_\mathrm{fn}$ & 3 & \\
        & & $a_\mathrm{sn}$ & 3 & \\
        & & $\kappa$ & 5e4 & \si{\pascal} \\

        \cline{2-5}
        & 
        & $K_\perp^\mathrm{epi}$ & 2e5 & \si{\pascal\per\metre} \\
        & Boundary & $K_\parallel^\mathrm{epi}$ & 2e4 & \si{\pascal\per\metre} \\
        & conditions & $C_\perp^\mathrm{epi}$ & 2e4 & \si{\pascal\second\per\metre} \\
        & & $C_\parallel^\mathrm{epi}$ & 2e3 & \si{\pascal\second\per\metre} \\

        \cline{2-5}
        & \multirow{1}{3.75cm}{In. conditions}
        & $p_0$ & 1333.2 & \si{\pascal} \\
        % \hline
    \end{tabular}
    \caption{Parameters used in the electromechanical model: electrophysiology (EP), active force generation (AFG) and solid mechanics (M). For the force generation model, we only report parameters that are different from the original setting described in \cite{regazzoni2020biophysically}. }
    \label{tab:params_ep}
\end{table}

\begin{table}
    \centering

    \begin{tabular}{l r S S}
        % \hline
        & \textbf{Parameter} & \multicolumn{2}{c}{\textbf{Value}} \\
        \hline\multirow{6}{*}{Systemic arteries}
        & $R_\text{AR}^\text{SYS}$ & 0.3750 & \siresistance \\
        & $C_\text{AR}^\text{SYS}$ & 2.048 & \sicapacitance\\
        & $L_\text{AR}^\text{SYS}$ &\num{2.7e-3} &  \siinductance\\
        & $R_\text{upstream}^\text{SYS}$ & 0.05 & \siresistance\\
        & $p_\text{AR}^\text{SYS}(0)$ & 80.0 & \si{\pascal}\\
        & $Q_\text{AR}^\text{SYS}(0)$ & 0.0 & \siflowrate \\
        
        \hline\multirow{5}{*}{Systemic veins}
        & $R_\text{VEN}^\text{SYS}$ & 0.26 & \siresistance \\
        & $C_\text{VEN}^\text{SYS}$ &60.0 & \sicapacitance\\
        & $L_\text{VEN}^\text{SYS}$ &\num{5e-4} & \siinductance\\
        & $p_\text{VEN}^\text{SYS}(0)$ & 30.9 & \si{\pascal}\\
        & $Q_\text{VEN}^\text{SYS}(0)$ & 0.0 & \si{\siflowrate}\\

        \hline\multirow{5}{*}{Pulmonary arteries}
        & $R_\text{AR}^\text{PUL}$ & 0.05& \siresistance \\
        & $C_\text{AR}^\text{PUL}$ & 10.0& \sicapacitance \\
        & $L_\text{AR}^\text{PUL}$ & \num{5e-4}& \siinductance \\
        & $p_\text{AR}^\text{PUL}(0)$ & 29.34 & \si{\pascal}\\
        & $Q_\text{AR}^\text{PUL}(0)$ &0.0 & \siflowrate \\

        \hline\multirow{5}{*}{Pulmonary veins}
        & $R_\text{VEN}^\text{PUL}$ & 0.025& \siresistance \\
        & $C_\text{VEN}^\text{PUL}$ & 38.4& \sicapacitance \\
        & $L_\text{VEN}^\text{PUL}$ &\num{2.083e-4} & \siinductance \\
        & $p_\text{VEN}^\text{PUL}(0)$ & 13.58& \si{\pascal} \\
        & $Q_\text{VEN}^\text{PUL}(0)$ & 0.0 & \siflowrate\\
        % \hline
            \end{tabular}

    \caption{Parameters of the circulation model for the ventricular electromechanical simulation: external circulation. The same parameters are employed for the 3D-0D CFD simulation.}
    \label{tab:params_circulation_I}
\end{table}

\begin{table}
    \centering
    \begin{tabular}{l r S S}
        % \hline
        & \textbf{Parameter} & \multicolumn{2}{c}{\textbf{Value}} \\
        \hline\multirow{6}{*}{Left atrium}
        & $E_\mathrm{A}$ & 0.07 & \sielastance \\
        & $E_\mathrm{B}$ & 0.09 & \sielastance \\
        & $t_\mathrm{C}$ & 0.80 & \\
        & $T_\mathrm{C}$ & 0.17 & \\
        & $T_\mathrm{R}$ & 0.17 & \\
        & $V_\mathrm{LA}(0)$ & 79.5& \si{\milli\litre} \\
        
        \hline\multirow{6}{*}{Right atrium}
        & $E_\mathrm{A}$ & 0.06 & \sielastance \\
        & $E_\mathrm{B}$ & 0.07 & \sielastance \\
        & $t_\mathrm{C}$ & 0.80 & \\
        & $T_\mathrm{C}$ & 0.17 & \\
        & $T_\mathrm{R}$ & 0.17 & \\
        & $V_\mathrm{RA}(0)$ & 64.17 & \si{\milli\litre} \\

        \hline\multirow{6}{*}{Right ventricle}
        & $E_\mathrm{A}$ & 0.55 & \sielastance \\
        & $E_\mathrm{B}$ & 0.05& \sielastance \\
        & $t_\mathrm{C}$ & 0.0 & \\
        & $T_\mathrm{C}$ & 0.34 & \\
        & $T_\mathrm{R}$ & 0.15 & \\
        & $V_\mathrm{RV}(0)$ & 148.9& \si{\milli\litre} \\
        
        \hline\multirow{2}{*}{Mitral valve}
        & $R_\mathrm{min}$ & 0.0164 & \si{\mmhg\second\per\milli\litre} \\
        & $R_\mathrm{max}$ & \num{75006.2} & \si{\mmhg\second\per\milli\litre} \\
        
        \hline\multirow{2}{*}{Aortic valve}
        & $R_\mathrm{min}$ &  0.0355 & \si{\mmhg\second\per\milli\litre} \\
        & $R_\mathrm{max}$ & \num{75006.2} & \si{\mmhg\second\per\milli\litre} \\
        
        \hline\multirow{2}{*}{Tricuspid valve}
        & $R_\mathrm{min}$ & 0.0075 & \si{\mmhg\second\per\milli\litre} \\
        & $R_\mathrm{max}$ & \num{75006.2} & \si{\mmhg\second\per\milli\litre} \\
        
        \hline\multirow{2}{*}{Pulmonary valve}
        & $R_\mathrm{min}$ & 0.0075 & \si{\mmhg\second\per\milli\litre} \\
        & $R_\mathrm{max}$ & \num{75006.2} & \si{\mmhg\second\per\milli\litre} \\

        % \hline
    \end{tabular}

    \caption{Parameters of the circulation model for the ventricular electromechanical simulation: cardiac circulation. Initial time of contraction $t_\mathrm{c}$, contraction duration $T_\mathrm{C}$ and relaxation duration $T_\mathrm{R}$ are relative to the heartbeat period. For the right atrium, right ventricle, tricuspid and pulmonary valves, the same parameters are employed for the 3D-0D CFD simulation.}
    \label{tab:params_circulation_II}
\end{table}

\FloatBarrier

%\nocite{*}% Show all bib entries - both cited and uncited; comment this line to view only cited bib entries;
%\bibliography{wileyNJD-AMA}%

%

\end{document}